# CAN ONE ESTIMATE THE CONDITIONAL DISTRIBUTION OF POST-MODEL-SELECTION ESTIMATORS?

By Hannes Leeb[1] and Benedikt M. Pötscher

*Yale University and University of Vienna*

We consider the problem of estimating the conditional distribution of a post-model-selection estimator where the conditioning is on the selected model. The notion of a post-model-selection estimator here refers to the combined procedure resulting from first selecting a model (e.g., by a model selection criterion such as AIC or by a hypothesis testing procedure) and then estimating the parameters in the selected model (e.g., by least-squares or maximum likelihood), all based on the same data set. We show that it is impossible to estimate this distribution with reasonable accuracy even asymptotically. In particular, we show that no estimator for this distribution can be uniformly consistent (not even locally). This follows as a corollary to (local) minimax lower bounds on the performance of estimators for this distribution. Similar impossibility results are also obtained for the conditional distribution of linear functions (e.g., predictors) of the post-model-selection estimator.

**1. Introduction and overview.** In many statistical applications a data-based model selection step precedes parameter estimation and inference. For example, the specification of the model (choice of functional form, choice of regressors, number of lags, etc.) is often based on the data. In contrast, the traditional theory of statistical inference is concerned with the properties of estimators and inference procedures under the central assumption of an a priori given model. That is, it is assumed that the model is known to the researcher prior to the statistical analysis, except for the value of the true

Received October 2003; revised November 2005.

[1]Supported by the Max Kade Foundation and by Austrian National Science Foundation (FWF) Grant P13868-MAT. A preliminary draft of the material in this paper was already written in 1999.

parameter vector. As a consequence, the actual statistical properties of estimators or inference procedures following a data-driven model selection step are not described by the traditional theory which assumes an a priori given model; in fact, they may differ substantially from the properties predicted by this theory; see, for example, [3, 4], [18], Section 3.3, or [21], Section 12. Ignoring the additional uncertainty originating from the data-driven model selection step and (inappropriately) applying traditional theory can hence result in very misleading conclusions.

Investigations into the distributional properties of post-model-selection estimators, that is, of estimators constructed after a data-driven model selection step, are relatively few and of recent vintage. Sen [23] obtained the unconditional large-sample limit distribution of a post-model-selection estimator in an i.i.d. maximum likelihood framework, when selection is between two competing nested models. In [18] the asymptotic properties of a class of post-model-selection estimators (based on a sequence of hypothesis tests) were studied in a rather general setting covering nonlinear models, dependent processes and more than two competing models. In that paper, the large-sample limit distribution of the post-model-selection estimator was derived, both unconditional as well as conditional on having chosen a correct model, not necessarily the minimal one. See also [20] for further discussion and a simulation study. The finite-sample distribution of a post-model-selection estimator, both unconditional and conditional on having chosen a particular (possibly incorrect) model, was derived in [12] in a normal linear regression framework; this paper also studied asymptotic approximations that are in a certain sense superior to the asymptotic distribution derived in [18]. The distributions of corresponding linear predictors constructed after model selection were studied in [10, 11]. Related work can also be found in [1, 5, 7, 8, 9, 15, 19, 24].

It transpires from the papers mentioned above that the finite-sample distributions (as well as the large-sample limit distributions) of post-model-selection estimators typically depend on unknown model parameters, often in a complicated fashion. For inference purposes, for example, for the construction of confidence sets, estimators for these distributions would be desirable. Consistent estimators for these distributions can typically be constructed quite easily, for example, by suitably replacing unknown parameters in the large-sample limit distributions by consistent estimators; see Section 2.2.1. However, the merits of such "plug-in" estimators in small samples are questionable: It is known that the convergence of the finite-sample distributions to their large-sample limits is typically not uniform with respect to the underlying parameters (see [12, 15] and Remark 4.11 in [14]), and there is no reason to believe that this nonuniformity will disappear when unknown parameters in the large-sample limit are replaced by estimators. This observation is the main motivation for the present paper to investigate in general the performance of estimators for the distribution of a post-model-selection estimator, where the estimators for the distribution are not necessarily "plug-in" estimators based on the limiting distribution. In particular, we ask whether estimators for the distribution function of post-model-selection estimators exist that do not suffer from the nonuniformity phenomenon mentioned above. As we show in this paper, the answer in general is "No." We also show that these negative results extend to the problem of estimating the distribution of linear functions (e.g., linear predictors) of post-model-selection estimators.

To fix ideas, consider for the moment the linear regression model

$$Y = V\chi + W\psi + u, \tag{1}$$

where $V$ and $W$, respectively, represent $n \times k$ and $n \times l$ nonstochastic regressor matrices ($k \geq 1, l \geq 1$), and the $n \times 1$ disturbance vector $u$ is normally distributed with mean zero and variance–covariance matrix $\sigma^2 I_n$. We also assume for the moment that $(V:W)'(V:W)/n$ converges to a nonsingular matrix as the sample size $n$ goes to infinity and that $\lim_{n\to\infty} V'W/n \neq 0$ (for a discussion of the case where this limit is zero, see Example 1 in Section 2.2.2). Now suppose that the vector $\chi$ represents the parameters of interest, while the parameter vector $\psi$ and the associated regressors in $W$ have been entered into the model only to avoid possible misspecification. Suppose further that the necessity to include $\psi$ or some of its components is then checked on the basis of the data, that is, a model selection procedure is used to determine which components of $\psi$ are to be retained in the model, the inclusion of $\chi$ not being disputed. The selected model is then used to obtain the final (post-model-selection) estimator $\tilde{\chi}$ for $\chi$. We are now interested in the conditional finite-sample distribution of $\tilde{\chi}$ (appropriately scaled and centered) given the outcome of the model selection step. (The reasons why we concentrate on the conditional rather than on the unconditional distribution are discussed below.) Denote this $k$-dimensional cumulative distribution function (c.d.f.) by $G_{n,\theta,\sigma}(t|\hat{M})$, where $\hat{M}$ stands for the selected model, that is, for the set of selected regressors. As indicated in the notation, this distribution function depends on the true parameters $\theta = (\chi', \psi')'$ and $\sigma$. For the sake of definiteness of discussion, assume for the moment that the model selection procedure used here is the particular "general-to-specific" procedure described at the beginning of Section 2; we comment on other model selection procedures, including Akaike's AIC and thresholding procedures, below.

As mentioned above, it is not difficult to construct a consistent estimator for $G_{n,\theta,\sigma}(t|\hat{M})$ for any given $t$, that is, an estimator $\hat{G}_n(t|\hat{M})$ satisfying

$$P_{n,\theta,\sigma}(|\hat{G}_n(t|\hat{M}) - G_{n,\theta,\sigma}(t|\hat{M})| > \delta) \stackrel{n\to\infty}{\longrightarrow} 0 \tag{2}$$

for each $\delta > 0$ and each $\theta$, $\sigma$; see Section 2.2.1. However, it follows from the results in Section 2.2.2 that *any* estimator satisfying (2), that is, *any consistent* estimator for $G_{n,\theta,\sigma}(t|\hat{M})$, necessarily also satisfies

$$(3) \qquad \liminf_{n\to\infty} \sup_{\|\theta\|<R} P_{n,\theta,\sigma}(|\hat{G}_n(t|\hat{M}) - G_{n,\theta,\sigma}(t|\hat{M})| > \delta) \geq c > 0$$

for suitable positive constants $c$, $R$ and $\delta$ (not depending on the estimator), with the lower bound $c$ often being quite large. That is, while the probability in (2) converges to zero for every given $\theta$ by consistency, relation (3) shows that it does not do so uniformly in $\theta$. It follows that $\hat{G}_n(t|\hat{M})$ can never be uniformly consistent (not even when restricting consideration to uniform consistency over all compact subsets of the parameter space). Hence, a large sample size does not guarantee a small estimation error with high probability when estimating the conditional distribution function of a post-model-selection estimator. In this sense, reliably assessing the precision of post-model-selection estimators is an intrinsically hard problem. Apart from (3), we also provide minimax lower bounds for arbitrary (not necessarily consistent) estimators of the conditional distribution function $G_{n,\theta,\sigma}(t|\hat{M})$. For example, we provide results that imply that

$$(4) \qquad \liminf_{n\to\infty} \inf_{\hat{G}_n(t|\hat{M})} \sup_{\|\theta\|<R} P_{n,\theta,\sigma}(|\hat{G}_n(t|\hat{M}) - G_{n,\theta,\sigma}(t|\hat{M})| > \delta) > 0$$

holds for suitable positive constants $R$ and $\delta$, where the infimum extends over *all* estimators of $G_{n,\theta,\sigma}(t|\hat{M})$. The results in Section 2.2.2 in fact show that the balls $\|\theta\| < R$ in (3) and (4) can be replaced by suitable balls (not necessarily centered at the origin) shrinking at the rate $n^{-1/2}$. This shows that the nonuniformity phenomenon described in (3)–(4) is a local, rather than a global, phenomenon. Moreover, relations (3)–(4) also hold with the unconditional probability $P_{n,\theta,\sigma}(\cdot)$ in (3)–(4) replaced by the conditional probability given model $M$ is selected, that is, given the event $\hat{M} = M$. In Section 2.2.2 we further show that the nonuniformity phenomenon expressed in (3) and (4) typically also arises when the parameter of interest is not $\chi$, but some other linear transformation of $\theta = (\chi', \psi')'$. As discussed in Remark 4.8, the results also hold for randomized estimators of the conditional distribution function $G_{n,\theta,\sigma}(t|\hat{M})$. Hence, no resampling procedure whatsoever can alleviate the problem. This explains the anecdotal evidence in the literature that resampling methods are often unsuccessful in approximating distributional properties of post-model-selection estimators (e.g., [4] or [6]).

The results outlined above are presented in Section 2.2.2 for the particular "general-to-specific" model selection procedure described at the beginning of Section 2. Analogous results for a large class of model selection procedures, including Akaike's AIC and thresholding procedures, are then given in Section 3 based on the results in Section 2.2.2. In fact, it transpires from the proofs that the nonuniformity phenomenon expressed in (3)–(4) is not specific to the model selection procedures discussed in Sections 2.2 and 3 of the present paper, but will occur for most (if not all) model selection procedures, including consistent ones; see Section 5.

In the present paper we focus on the *conditional* distribution of the post-model-selection estimator. Given that the outcome of the model selector has been observed, it may be argued that the relevant sample space for assessing variability of the parameter estimator is then not given by the entire original sample space, but rather by the subset that gave rise to the observed outcome of the model selector; see the literature on conditional inference ([22] and [17], page 421). If one does not adhere to such a conditionality principle the unconditional distribution of the post-model-selection estimator is of interest. For this case, similar results can be obtained and are reported in [13].

The plan of the paper is as follows: Post-model-selection estimators based on a "general-to-specific" model selection procedure are the subject of Section 2. After introducing the basic framework and some notation, such as the family of models $M_p$ from which the "general-to-specific" model selection procedure $\hat{p}$ selects as well as the post-model-selection estimator $\tilde{\theta}$, the conditional c.d.f. $G_{n,\theta,\sigma}(t|p)$ of (a linear function of) the post-model-selection estimator $\tilde{\theta}$ given that $\hat{p}$ selects model $M_p$ is introduced and discussed in Section 2.1. Consistent estimation of $G_{n,\theta,\sigma}(t|p)$ and of $G_{n,\theta,\sigma}(t|\hat{p})$ (i.e., of the c.d.f. conditional on the actual outcome of the model selection procedure) is discussed in Section 2.2.1. The main results of the paper are contained in Sections 2.2.2 and 3: In Section 2.2.2 we provide a detailed analysis of the nonuniformity phenomenon encountered in (3)–(4). In Section 3 the "impossibility" result from Section 2.2.2 is extended to a large class of model selection procedures, including Akaike's AIC, and to selection from a nonnested collection of models. Further theoretical results on which the proofs are based are given in Section 4 and conclusions are drawn in Section 5. All proofs, as well as some auxiliary results, are collected into appendices. Finally, a word on notation: The Euclidean norm is denoted by $\|\cdot\|$, and $\lambda_{\max}(E)$ denotes the largest eigenvalue of a symmetric matrix $E$. A prime denotes transposition of a matrix. For vectors $x$ and $y$, the relation $x \leq y$ ($x < y$, resp.) denotes $x_i \leq y_i$ ($x_i < y_i$, resp.) for all $i$. As usual, $\Phi$ denotes the standard normal distribution function.

**2. Results for post-model-selection estimators based on a "general-to-specific" model selection procedure.** Consider the linear regression model

(5) $$Y = X\theta + u,$$

where $X$ is a nonstochastic $n \times P$ matrix with $\mathrm{rank}(X) = P$ and $u \sim N(0, \sigma^2 I_n)$, $\sigma^2 > 0$. Here $n$ denotes the sample size and we assume $n > P \geq 1$. In addition, we assume that $Q = \lim_{n \to \infty} X'X/n$ exists and is nonsingular. In this section we shall—similarly as in [18]—consider model selection from the collection of nested models $M_\mathcal{O} \subseteq M_{\mathcal{O}+1} \subseteq \cdots \subseteq M_P$, where $\mathcal{O}$ is specified by the user, and where, for $0 \leq p \leq P$, the model $M_p$ is given by

$$M_p = \{(\theta_1, \ldots, \theta_P)' \in \mathbf{R}^P : \theta_{p+1} = \cdots = \theta_P = 0\}.$$

(In Section 3 below general nonnested families of models will also be considered.) Clearly, the model $M_p$ corresponds to the situation where only the first $p$ regressors in (5) are included. For the most parsimonious model under consideration, that is, for $M_\mathcal{O}$, we assume that $\mathcal{O}$ satisfies $0 \leq \mathcal{O} < P$; if $\mathcal{O} > 0$, this model contains as free parameters only those components of the parameter vector $\theta$ that are not subject to model selection. [In the notation used in connection with (1) we then have $\chi = (\theta_1, \ldots, \theta_\mathcal{O})'$ and $\psi = (\theta_{\mathcal{O}+1}, \ldots, \theta_P)'$.] Furthermore, note that $M_0 = \{(0, \ldots, 0)'\}$ and that $M_P = \mathbf{R}^P$. We call $M_p$ the regression model of order $p$.

The following notation will prove useful. For matrices $B$ and $C$ of the same row-dimension, the column-wise concatenation of $B$ and $C$ is denoted by $(B:C)$. If $D$ is an $m \times P$ matrix, let $D[p]$ denote the $m \times p$ matrix consisting of the first $p$ columns of $D$. Similarly, let $D[\neg p]$ denote the $m \times (P-p)$ matrix consisting of the last $P - p$ columns of $D$. If $x$ is a $P \times 1$ vector, we write in abuse of notation $x[p]$ and $x[\neg p]$ for $(x'[p])'$ and $(x'[\neg p])'$, respectively. (We shall use the above notation also in the "boundary" cases $p = 0$ and $p = P$. It will always be clear from the context how expressions containing symbols such as $D[0]$, $D[\neg P]$, $x[0]$ or $x[\neg P]$ are to be interpreted.) As usual, the $i$th component of a vector $x$ is denoted by $x_i$, and the entry in the $i$th row and $j$th column of a matrix $B$ is denoted by $B_{i,j}$.

The restricted least-squares estimator of $\theta$ under the restriction $\theta[\neg p] = 0$, that is, under $\theta_{p+1} = \cdots = \theta_P = 0$, will be denoted by $\tilde{\theta}(p)$, $0 \leq p \leq P$ (in case $p = P$ the restriction is void). Note that $\tilde{\theta}(p)$ is given by the $P \times 1$ vector

$$\tilde{\theta}(p) = \begin{pmatrix} (X[p]'X[p])^{-1}X[p]'Y \\ (0, \ldots, 0)' \end{pmatrix},$$

where the expressions $\tilde{\theta}(0)$ and $\tilde{\theta}(P)$, respectively, are to be interpreted as the zero-vector in $\mathbf{R}^P$ and as the unrestricted least-squares estimator of $\theta$. Given a parameter vector $\theta$ in $\mathbf{R}^P$, the order of $\theta$ (relative to the nested sequence of models $M_p$) is defined as

$$p_0(\theta) = \min\{p : 0 \leq p \leq P, \theta \in M_p\}.$$

Hence, if $\theta$ is the true parameter vector, a model $M_p$ is a correct model if and only if $p \geq p_0(\theta)$. We stress that $p_0(\theta)$ is a property of a *single parameter*, and needs to be distinguished from the notion of the order of the model $M_p$ introduced earlier, which is a property of the *set of parameters* $M_p$.

A model selection procedure is now nothing else than a data-driven (measurable) rule $\hat{p}$ that selects a value from $\{\mathcal{O}, \ldots, P\}$ and thus selects a model from the list of candidate models $M_\mathcal{O}, \ldots, M_P$. In this section we shall consider as an important leading case a "general-to-specific" model selection procedure based on a sequence of hypothesis tests. (Results for a larger class of model selection procedures, including Akaike's AIC, are provided in Section 3.) This procedure is given as follows: The sequence of hypotheses $H_0^p : p_0(\theta) < p$ is tested against the alternatives $H_1^p : p_0(\theta) = p$ in decreasing order starting at $p = P$. If, for some $p > \mathcal{O}$, $H_0^p$ is the first hypothesis in the process that is rejected, we set $\hat{p} = p$. If no rejection occurs until even $H_0^{\mathcal{O}+1}$ is not rejected, we set $\hat{p} = \mathcal{O}$. Each hypothesis in this sequence is tested by a kind of $t$-test where the error variance is always estimated from the overall model (but see the discussion following Theorem 3.1 in Section 3 below for other choices of estimators of the error variance). More formally, we have

$$\hat{p} = \max\{p : |T_p| \geq c_p, 0 \leq p \leq P\}, \tag{6}$$

with $c_\mathcal{O} = 0$ in order to ensure a well-defined $\hat{p}$ in the range $\{\mathcal{O}, \mathcal{O}+1, \ldots, P\}$. For $\mathcal{O} < p \leq P$, the critical values $c_p$ satisfy $0 < c_p < \infty$ and are independent of sample size (but see also Remark 4.7). The test-statistics are given by

$$T_p = \frac{\sqrt{n}\tilde{\theta}_p(p)}{\hat{\sigma}\xi_{n,p}} \qquad (0 < p \leq P),$$

with the convention that $T_0 = 0$. Furthermore,

$$\xi_{n,p} = \left(\left[\left(\frac{X[p]'X[p]}{n}\right)^{-1}\right]_{p,p}\right)^{1/2} \qquad (0 < p \leq P)$$

denotes the nonnegative square root of the $p$th diagonal element of the matrix indicated, and $\hat{\sigma}^2$ is given by

$$\hat{\sigma}^2 = (n-P)^{-1}(Y - X\tilde{\theta}(P))'(Y - X\tilde{\theta}(P)).$$

Note that, under the hypothesis $H_0^p$, the statistic $T_p$ is $t$-distributed with $n - P$ degrees of freedom for $0 < p \leq P$. It is also easy to see that the so-defined model selection procedure $\hat{p}$ is conservative: The probability of selecting an incorrect model, that is, the probability of the event $\{\hat{p} < p_0(\theta)\}$, converges to zero as the sample size increases. In contrast, the probability of selecting a correct (but possibly overparameterized) model, that is, the probability of the event $\{\hat{p} = p\}$ for $p$ satisfying $\max\{p_0(\theta), \mathcal{O}\} \leq p \leq P$, converges to a positive limit; see, for example, Proposition 5.4 and equation (5.7) in [11].

The post-model-selection estimator $\tilde{\theta}$ can now be defined as follows: On the event $\hat{p} = p$, $\tilde{\theta}$ is given by the restricted least-squares estimator $\tilde{\theta}(p)$, that is,

$$\tilde{\theta} = \sum_{p=\mathcal{O}}^{P} \tilde{\theta}(p)\, \mathbf{1}(\hat{p} = p), \tag{7}$$

where $\mathbf{1}(\cdot)$ denotes the indicator function of the event shown in the argument.

2.1. *The conditional finite-sample distribution of the post-model-selection estimator.* We now introduce the distribution function of a linear transformation of $\tilde{\theta}$, conditional on the event $\hat{p} = p$, and summarize some of its properties that will be needed in the subsequent development. To this end, let $A$ be a nonstochastic $k \times P$ matrix of rank $k$, $1 \leq k \leq P$. For $\mathcal{O} \leq p \leq P$, we consider the conditional c.d.f.

$$G_{n,\theta,\sigma}(t|p) = P_{n,\theta,\sigma}(\sqrt{n} A(\tilde{\theta} - \theta) \leq t | \hat{p} = p) \qquad (t \in \mathbf{R}^k). \tag{8}$$

Here $P_{n,\theta,\sigma}(\cdot)$ denotes the probability measure corresponding to a sample of size $n$ from (5), and $P_{n,\theta,\sigma}(\cdot|\hat{p} = p)$ denotes the associated conditional probability measure (the conditioning event always having positive probability; cf. (3.8)–(3.9) in [11] and the attending discussion). Note that, on the event $\hat{p} = p$, the expression $A(\tilde{\theta} - \theta)$ equals $A(\tilde{\theta}(p) - \theta)$ in view of (7).

Depending on the choice of the matrix $A$, several important scenarios are covered by (8): The conditional c.d.f. of $\sqrt{n}(\tilde{\theta} - \theta)$ is obtained by setting $A$ equal to the $P \times P$ identity matrix $I_P$. The conditional c.d.f. of the components of $\sqrt{n}(\tilde{\theta} - \theta)$ that are not restricted to zero in the selected model $M_p$, $p > 0$, is obtained by setting $A$ to the $p \times P$ matrix $(I_p : 0)$. In case $\mathcal{O} > 0$, the conditional c.d.f. of those components of $\sqrt{n}(\tilde{\theta} - \theta)$ which correspond to the parameter of interest $\chi$ in (1) can be studied by setting $A$ to the $\mathcal{O} \times P$ matrix $(I_{\mathcal{O}} : 0)$, as we then have $A\theta = (\theta_1, \ldots, \theta_{\mathcal{O}})' = \chi$. Finally, if $A \neq 0$ is a $1 \times P$ vector, we obtain the conditional distribution of a linear predictor based on the post-model-selection estimator. See the examples at the ends of Section 2.2.2 and Section 4.1 for more discussion.

The c.d.f. $G_{n,\theta,\sigma}(t|p)$ and its properties have been analyzed in detail in [12] and [10]. To be able to access these results, we need some further notation. The expected value of the restricted least-squares estimator $\tilde{\theta}(p)$ will be denoted by $\eta_n(p)$ and is given by the $P \times 1$ vector

$$\eta_n(p) = \begin{pmatrix} \theta[p] + (X[p]'X[p])^{-1}X[p]'X[\neg p]\theta[\neg p] \\ (0, \ldots, 0)' \end{pmatrix}, \tag{9}$$

with the conventions that $\eta_n(0) = (0, \ldots, 0)' \in \mathbf{R}^P$ and that $\eta_n(P) = \theta$. Furthermore, let $\Phi_{n,p}(\cdot)$ denote the c.d.f. of $\sqrt{n} A(\tilde{\theta}(p) - \eta_n(p))$, that is, the c.d.f. of $\sqrt{n}A$ times the restricted least-squares estimator based on model $M_p$ centered at its mean. Hence, $\Phi_{n,p}(\cdot)$ is the c.d.f. of a $k$-variate Gaussian random vector with mean zero and variance–covariance matrix $\sigma^2 A[p](X[p]'X[p]/n)^{-1}A[p]'$ in case $p > 0$, and it is the c.d.f. of point-mass at zero in $\mathbf{R}^k$ in case $p = 0$. If $p > 0$ and if the matrix $A[p]$ has full row rank $k$, then $\Phi_{n,p}(\cdot)$ has a density with respect to Lebesgue measure, and we shall denote this density by $\phi_{n,p}(\cdot)$. We note that $\eta_n(p)$ depends on $\theta$ and that $\Phi_{n,p}(\cdot)$ depends on $\sigma$ (in case $p > 0$), although these dependencies are not shown explicitly in the notation.

For $p > 0$, we introduce

$$b_{n,p} = C_n^{(p)'}(A[p](X[p]'X[p]/n)^{-1}A[p]')^- \tag{10}$$

and

$$\zeta_{n,p}^2 = \xi_{n,p}^2 - C_n^{(p)'}(A[p](X[p]'X[p]/n)^{-1}A[p]')^- C_n^{(p)}, \tag{11}$$

with $C_n^{(p)} = A[p](X[p]'X[p]/n)^{-1}e_p$, where $e_p$ denotes the $p$th standard basis vector in $\mathbf{R}^p$, and $B^-$ denotes a generalized inverse of the matrix $B$. (Observe that $\zeta_{n,p}^2$ is invariant under the choice of the generalized inverse. The same is not necessarily true for $b_{n,p}$, but is true for $b_{n,p}z$ for all $z$ in the column-space of $A[p]$. Also note that (13) below depends on $b_{n,p}$ only through $b_{n,p}z$ with $z$ in the column-space of $A[p]$.) We observe that the vector of covariances between $A\tilde{\theta}(p)$ and $\tilde{\theta}_p(p)$ is precisely given by $\sigma^2 n^{-1} C_n^{(p)}$ (and hence does *not* depend on $\theta$). Furthermore, observe that $A\tilde{\theta}(p)$ and $\tilde{\theta}_p(p)$ are uncorrelated if and only if $\zeta_{n,p}^2 = \xi_{n,p}^2$ if and only if $b_{n,p}z = 0$ for all $z$ in the column-space of $A[p]$; see Lemma A.2 in [10].

Finally, for a univariate Gaussian random variable $\mathfrak{N}$ with zero mean and variance $s^2$, $s \geq 0$, we write $\Delta_s(a,b)$ for $P(|\mathfrak{N} - a| < b)$, $a \in \mathbf{R} \cup \{-\infty, \infty\}$, $b \in \mathbf{R}$. Note that $\Delta_s(\cdot, \cdot)$ is symmetric around zero in its first argument, and that $\Delta_s(-\infty, b) = \Delta_s(\infty, b) = 0$ holds. In case $s = 0$, $\mathfrak{N}$ is to be interpreted as being equal to zero, hence, $a \mapsto \Delta_0(a,b)$ reduces to the indicator function of the interval $(-b, b)$.

We are now in a position to present the explicit formulae for $G_{n,\theta,\sigma}(t|p)$ derived in [10]. In case $p = \mathcal{O}$ we have

$$G_{n,\theta,\sigma}(t|\mathcal{O}) = \Phi_{n,\mathcal{O}}(t - \sqrt{n}A(\eta_n(\mathcal{O}) - \theta)), \tag{12}$$

that is, the c.d.f. of (a linear function of) the post-model-selection estimator $\tilde{\theta}$ conditional on $\hat{p} = \mathcal{O}$ coincides with the c.d.f. of (this linear function of) the restricted least-squares estimator $\tilde{\theta}(\mathcal{O})$. However, in case $p > \mathcal{O}$ we have

$$G_{n,\theta,\sigma}(t|p) = \int_{z \leq t - \sqrt{n}A(\eta_n(p) - \theta)} m_{n,p,\theta,\sigma}(z) \Phi_{n,p}(dz). \tag{13}$$

In the above display, $\Phi_{n,p}(dz)$ denotes integration with respect to the measure induced by the normal c.d.f. $\Phi_{n,p}(\cdot)$ on $\mathbf{R}^k$ and the integrand $m_{n,p,\theta,\sigma}(z)$

is given by

$$m_{n,p,\theta,\sigma}(z) = \left[ \int_0^\infty (1 - \Delta_{\sigma\zeta_{n,p}}(\sqrt{n}\eta_{n,p}(p) + b_{n,p}z, sc_p\sigma\xi_{n,p})) \right.$$
$$\left. \times \prod_{q=p+1}^P \Delta_{\sigma\xi_{n,q}}(\sqrt{n}\eta_{n,q}(q), sc_q\sigma\xi_{n,q})h(s)\,ds \right] \quad (14)$$
$$\times [P_{n,\theta,\sigma}(\hat{p}=p)]^{-1},$$

where $\zeta_{n,p}$ is the nonnegative root of $\zeta_{n,p}^2$ and the model selection probability $P_{n,\theta,\sigma}(\hat{p}=p)$ is given by

$$P_{n,\theta,\sigma}(\hat{p}=p) = \left[ \int_0^\infty (1 - \Delta_{\sigma\xi_{n,p}}(\sqrt{n}\eta_{n,p}(p), sc_p\sigma\xi_{n,p})) \right.$$
$$\left. \times \prod_{q=p+1}^P \Delta_{\sigma\xi_{n,q}}(\sqrt{n}\eta_{n,q}(q), sc_q\sigma\xi_{n,q})h(s)\,ds \right]. \quad (15)$$

In the two displays above, $h$ denotes the density of $\hat{\sigma}/\sigma$, that is, $h$ is the density of $(n-P)^{-1/2}$ times the square-root of a chi-square distributed random variable with $n - P$ degrees of freedom. The conditional finite-sample distribution of the post-model-selection estimator given in (13) is not normal; for example, it can be bimodal; see Figure 2 in [11]. An exception where (13) is normal is the case where $C_n^{(p)} = 0$, that is, when $A\tilde{\theta}(p)$ and $\tilde{\theta}_p(p)$ are uncorrelated; see [10], Section 3.3, for more discussion. On the other extreme, namely, if $A\tilde{\theta}(p)$ and $\tilde{\theta}_p(p)$ are perfectly correlated in the sense that $\zeta_{n,p} = 0$ holds, the function $\Delta_{\sigma\zeta_{n,p}}$ appearing in (14) reduces to an indicator function. This is, for example, the case if $A = I_P$ or if $A = (I_p:0)$.

2.2. *Estimators of the conditional finite-sample distribution.* For the purpose of inference after model selection, the conditional finite-sample distribution of the post-model-selection-estimator is an object of particular interest. As we have seen, it depends on unknown parameters in a complicated manner, and, hence, one will have to be satisfied with estimators of this c.d.f. The object we would primarily like to estimate is

$$G_{n,\theta,\sigma}(t|\hat{p}) = \sum_{p=\mathcal{O}}^P G_{n,\theta,\sigma}(t|p)\mathbf{1}(\hat{p}=p),$$

that is, the conditional c.d.f. after the model selection procedure has returned the model order $\hat{p}$. As we shall see in Section 2.2.1, it is not difficult to construct consistent estimators for $G_{n,\theta,\sigma}(t|\hat{p})$. We note that in considering consistency of an estimator of $G_{n,\theta,\sigma}(t|\hat{p})$ one is evaluating the performance

of such an estimator in an unconditional manner, namely, over the entire sample space. One can also take a conditional view in such an evaluation and ask if the given estimator of $G_{n,\theta,\sigma}(t|\hat{p})$ is "consistent conditionally on the outcome $\hat{p} = p$," at least for those parameter values $\theta$ that lead to a positive limit of $P_{n,\theta,\sigma}(\hat{p} = p)$, which are precisely all $\theta \in M_p$ as shown in Proposition A.2 in Appendix A. Of course, this reduces then to the question of (conditional) consistency of estimators of $G_{n,\theta,\sigma}(t|p)$ also discussed in Section 2.2.1 below.

Despite the consistency results in Section 2.2.1, we shall find in Section 2.2.2 that *any* estimator of $G_{n,\theta,\sigma}(t|\hat{p})$ typically performs unsatisfactorily, in that the estimation error cannot become small uniformly over (subsets of) the parameter space even as sample size goes to infinity. In particular, no uniformly consistent estimators exist, not even locally. These results rest on parallel results for the estimation of $G_{n,\theta,\sigma}(t|p)$ with fixed $p$ which are collected in Section 4.1 below.

2.2.1. *Consistent estimators.* We construct consistent estimators for $G_{n,\theta,\sigma}(t|\hat{p})$ and $G_{n,\theta,\sigma}(t|p)$ (consistent over $M_p$ in the latter case) by commencing from the asymptotic distribution. The large-sample limit of $G_{n,\theta,\sigma}(t|p)$ for $\theta \in M_p$ is given by $G_{\infty,\theta,\sigma}(t|p) = \Phi_{\infty,p}(t)$ in case $p = \max\{p_0(\theta), \mathcal{O}\}$, and by

$$(16) \quad G_{\infty,\theta,\sigma}(t|p) = \int_{\substack{z \in \mathbf{R}^k \\ z \leq t}} \frac{1 - \Delta_{\sigma\zeta_{\infty,p}}(b_{\infty,p}z, c_p\sigma\xi_{\infty,p})}{1 - \Delta_{\sigma\xi_{\infty,p}}(0, c_p\sigma\xi_{\infty,p})} \Phi_{\infty,p}(dz)$$

in case $p > \max\{p_0(\theta), \mathcal{O}\}$. This follows from Proposition A.1 in Appendix A with $\gamma = 0$ and $\sigma^{(n)} = \sigma$. Here $\Phi_{\infty,p}$ is the c.d.f. of a $k$-variate Gaussian random vector with mean zero and variance–covariance matrix $\sigma^2 A[p]Q[p:p]^{-1}A[p]'$, $0 < p \leq P$, where $Q[p:p]$ represents the leading diagonal $p \times p$ submatrix of $Q$. Also, let $\Phi_{\infty,0}(\cdot)$ denote the c.d.f. of point-mass at zero in $\mathbf{R}^k$. Note that $G_{\infty,\theta,\sigma}(t|p)$, for $p > \mathcal{O}$, depends on $\theta$ as it follows two different formulas depending on whether $\theta \in M_p \setminus M_{p-1}$ or $\theta \in M_{p-1}$. Let $\hat{\Phi}_{n,p}(\cdot)$ denote the c.d.f. of a $k$-variate Gaussian random vector with mean zero and variance–covariance matrix $\hat{\sigma}^2 A[p](X[p]'X[p]/n)^{-1}A[p]'$, $0 < p \leq P$; we also adopt the convention that $\hat{\Phi}_{n,0}(\cdot)$ denotes the c.d.f. of point-mass at zero in $\mathbf{R}^k$. [We use the same convention for $\hat{\Phi}_{n,p}(\cdot)$ in case $\hat{\sigma} = 0$, which is a probability zero event.] For given $p$, $\mathcal{O} \leq p \leq P$, an estimator $\check{G}_n(t|p)$ for $G_{n,\theta,\sigma}(t|p)$ is now defined as follows: For $p = \mathcal{O}$, we set $\check{G}_n(t|\mathcal{O}) = \hat{\Phi}_{n,\mathcal{O}}(t)$. For $p > \mathcal{O}$, we first employ an auxiliary procedure that consistently decides between $p_0(\theta) = p$ and $p_0(\theta) < p$, that is, between $\theta \in M_p \setminus M_{p-1}$ and $\theta \in M_{p-1}$, for every $\theta \in M_p$. [E.g., the procedure that decides for $p_0(\theta) = p$ whenever $|T_p| > s_{n,p}$ and for $p_0(\theta) < p$ otherwise, with $s_{n,p}$ satisfying $s_{n,p} \to \infty$, $s_{n,p} = o(n^{1/2})$ for $n \to \infty$ can be used. Alternatively, a consistent model selection procedure such as BIC could be employed to select between $M_{p-1}$ and $M_p\backslash M_{p-1}$.] If the procedure decides for $p_0(\theta) = p$, we set $\check{G}_n(t|p) = \hat{\Phi}_{n,p}(t)$; otherwise we set $\check{G}_n(t|p)$ equal to the expression in (16) with $\hat{\sigma}$, $b_{n,p}$, $\zeta_{n,p}$, $\xi_{n,p}$ and $\hat{\Phi}_{n,p}(\cdot)$ replacing $\sigma$, $b_{\infty,p}$, $\zeta_{\infty,p}$, $\xi_{\infty,p}$ and $\Phi_{\infty,p}(\cdot)$, respectively. A little reflection shows that $\check{G}_n(t|p)$ is again a c.d.f. (This is trivial if $\hat{\sigma} = 0$, and follows for $\hat{\sigma} > 0$ from the observation that then $\check{G}_n(t|p)$ is either a normal c.d.f. or coincides with the conditional c.d.f. $G^*_{n,\theta,\sigma}(t|p)$ given in (13) of [10] with $\sigma$ replaced by $\hat{\sigma}$.) This gives an estimator $\check{G}_n(t|p)$ of $G_{n,\theta,\sigma}(t|p)$; as an estimator of $G_{n,\theta,\sigma}(t|\hat{p})$, we shall use $\check{G}_n(t|\hat{p}) = \sum_{p=\mathcal{O}}^{P} \check{G}_n(t|p) \mathbf{1}(\hat{p} = p)$. We have the following consistency results.

PROPOSITION 2.1. *Let $p$ satisfy $\mathcal{O} \leq p \leq P$. Then the estimator $\check{G}_n(t|p)$ is consistent (in the total variation distance) for $G_{n,\theta,\sigma}(t|p)$ and $G_{\infty,\theta,\sigma}(t|p)$ over the subset $M_p$ (and over $0 < \sigma < \infty$). That is, for every $\delta > 0$,*

(17) $$P_{n,\theta,\sigma}(\|\check{G}_n(\cdot|p) - G_{n,\theta,\sigma}(\cdot|p)\|_{\mathrm{TV}} > \delta) \stackrel{n \to \infty}{\longrightarrow} 0,$$

(18) $$P_{n,\theta,\sigma}(\|\check{G}_n(\cdot|p) - G_{\infty,\theta,\sigma}(\cdot|p)\|_{\mathrm{TV}} > \delta) \stackrel{n \to \infty}{\longrightarrow} 0$$

*for all $\theta \in M_p$ and all $\sigma > 0$. The results* (17) *and* (18) *also hold with $P_{n,\theta,\sigma}(\cdot|\hat{p} = p)$ replacing $P_{n,\theta,\sigma}(\cdot)$.*

COROLLARY 2.2. *The estimator $\check{G}_n(t|\hat{p})$ is consistent (in the total variation distance) for $G_{n,\theta,\sigma}(t|\hat{p})$ over the entire parameter space, that is, for every $\delta > 0$,*

$$P_{n,\theta,\sigma}(\|\check{G}_n(\cdot|\hat{p}) - G_{n,\theta,\sigma}(\cdot|\hat{p})\|_{\mathrm{TV}} > \delta) \stackrel{n \to \infty}{\longrightarrow} 0$$

*for all $\theta \in \mathbf{R}^P$ and $\sigma > 0$.*

While the estimators constructed above are consistent, they can be expected to perform poorly in finite samples when the true $\theta$ belongs to $M_p \backslash M_{p-1}$ but is "close" to $M_{p-1}$, since the auxiliary decision procedure (although being consistent) will then have difficulties making the correct decision in finite samples, and since $G_{n,\theta,\sigma}(\cdot|p)$ typically does not converge uniformly with respect to $\theta \in M_p \backslash M_{p-1}$ "close" to $M_{p-1}$ (cf. [12, 15] and Remark 4.11 in [14]). In the next section we show that this poor performance is not particular to the estimators constructed above, but is a genuine feature of the estimation problem under consideration.

2.2.2. *Performance limits and an impossibility result.* We now provide lower bounds for the performance of estimators of the conditional c.d.f. $G_{n,\theta,\sigma}(t|\hat{p})$ of the post-model-selection estimator $A\tilde{\theta}$; that is, we give lower bounds on the worst-case probability that the estimation error exceeds a certain threshold. These lower bounds are often quite large; furthermore, they remain lower bounds even if one restricts attention only to certain subsets of the parameter space that shrink at the rate $n^{-1/2}$. In this sense the "impossibility" results are of a local nature. In particular, the lower bounds imply that no uniformly consistent estimator of the conditional c.d.f. $G_{n,\theta,\sigma}(t|\hat{p})$ exists, not even locally. Similar results under a conditional evaluation of the estimation error are given in Section 4.1 and form the theoretical backbone for the results in the present section. We note already here that the lower bounds obtained in Section 4.1 are as large as 1 or 1/2, depending on the particular situation considered.

In the following, the asymptotic "correlation" between $A\tilde{\theta}(p)$ and $\tilde{\theta}_p(p)$ as measured by $C_\infty^{(p)} = \lim_{n\to\infty} C_n^{(p)}$ will play an important rôle. Note that $C_\infty^{(p)}$ equals $A[p]Q[p:p]^{-1}e_p$, and hence, does *not* depend on the unknown parameters $\theta$ or $\sigma$. In the important special case discussed in the Introduction [cf. (1)], the matrix $A$ equals the $\mathcal{O} \times P$ matrix $(I_\mathcal{O}:0)$, and the condition $C_\infty^{(p)} \neq 0$ reduces to the condition that the regressor corresponding to the $p$th column of $(V:W)$ is asymptotically correlated with at least one of the regressors corresponding to the columns of $V$. See Example 1 below for more discussion.

In the result to follow we shall consider performance limits for estimators of $G_{n,\theta,\sigma}(t|\hat{p})$ at a *fixed* value of the argument $t$. An estimator of $G_{n,\theta,\sigma}(t|\hat{p})$ is now nothing else than a real-valued random variable $\Gamma_n = \Gamma_n(Y,X)$. For mnemonic reasons, we shall, however, use the symbol $\hat{G}_n(t|\hat{p})$ instead of $\Gamma_n$ to denote an arbitrary estimator of $G_{n,\theta,\sigma}(t|\hat{p})$. This notation should not be taken as implying that the estimator is obtained by evaluating an estimated c.d.f. at the argument $t$, or that it is a priori constrained to lie between zero and one. We shall use this notational convention mutatis mutandis also in subsequent sections. Regarding the nonuniformity phenomenon, we then have a dichotomy which is described in the following two results.

THEOREM 2.3. *Suppose that $A\tilde{\theta}(q)$ and $\tilde{\theta}_q(q)$ are asymptotically correlated, that is, $C_\infty^{(q)} \neq 0$, for some $q$ satisfying $\mathcal{O} < q \leq P$, and let $q^*$ denote the largest $q$ with this property. Then the following holds for each $\theta \in M_{q^*-1}$, $0 < \sigma < \infty$, and each $t \in \mathbf{R}^k$: There exist $\delta_0 > 0$ and $0 < \rho_0 < \infty$ such that any estimator $\hat{G}_n(t|\hat{p})$ of $G_{n,\theta,\sigma}(t|\hat{p})$ satisfying*

$$(19) \qquad P_{n,\theta,\sigma}(|\hat{G}_n(t|\hat{p}) - G_{n,\theta,\sigma}(t|\hat{p})| > \delta) \stackrel{n\to\infty}{\longrightarrow} 0$$

*for each $\delta > 0$ (in particular, every estimator that is consistent) also satisfies*

$$\liminf_{n\to\infty} \sup_{\substack{\vartheta \in M_{q^*} \\ \|\vartheta-\theta\| < \rho_0/\sqrt{n}}} P_{n,\vartheta,\sigma}(|\hat{G}_n(t|\hat{p}) - G_{n,\vartheta,\sigma}(t|\hat{p})| > \delta_0)$$

(20)
$$\geq 2(1-\Phi(c_{q^*}))\prod_{q=q^*+1}^{P}(2\Phi(c_q)-1) > 0.$$

*The constants $\delta_0$ and $\rho_0$ may be chosen in such a way that they depend only on $t, Q, A, \sigma$ and the critical value $c_{q^*}$. Moreover,*

(21)   $$\liminf_{n\to\infty}\inf_{\hat{G}_n(t|\hat{p})}\sup_{\substack{\vartheta\in M_{q^*}\\ \|\vartheta-\theta\|<\rho_0/\sqrt{n}}} P_{n,\vartheta,\sigma}(|\hat{G}_n(t|\hat{p})-G_{n,\vartheta,\sigma}(t|\hat{p})|>\delta_0) > 0$$

*and*

(22)
$$\sup_{\delta>0}\liminf_{n\to\infty}\inf_{\hat{G}_n(t|\hat{p})}\sup_{\substack{\vartheta\in M_{q^*}\\ \|\vartheta-\theta\|<\rho_0/\sqrt{n}}} P_{n,\vartheta,\sigma}(|\hat{G}_n(t|\hat{p})-G_{n,\vartheta,\sigma}(t|\hat{p})|>\delta)$$
$$\geq (1-\Phi(c_{q^*}))\prod_{q=q^*+1}^{P}(2\Phi(c_q)-1) > 0$$

*hold, where the infima in* (21) *and* (22) *extend over* all *estimators $\hat{G}_n(t|\hat{p})$ of $G_{n,\theta,\sigma}(t|\hat{p})$.* [*The lower bound in* (20) *is nothing else than* $\lim_{n\to\infty} P_{n,\theta,\sigma}(\hat{p}=q^*)$.]

PROPOSITION 2.4. *Suppose that $A\tilde{\theta}(q)$ and $\tilde{\theta}_q(q)$ are asymptotically uncorrelated, that is, $C_{\infty}^{(q)}=0$, for all $q$ satisfying $\mathcal{O}<q\leq P$. Then*

(23)   $$\sup_{\theta\in\mathbf{R}^P}\sup_{\substack{\sigma\in\mathbf{R}\\ \sigma_*\leq\sigma\leq\sigma^*}} P_{n,\theta,\sigma}(\|\hat{\tilde{\Phi}}_{n,P}(\cdot)-G_{n,\theta,\sigma}(\cdot|\hat{p})\|_{\mathrm{TV}}>\delta) \stackrel{n\to\infty}{\longrightarrow} 0$$

*holds for each $\delta>0$, and for any constants $\sigma_*$ and $\sigma^*$ satisfying $0<\sigma_*\leq\sigma^*<\infty$.*

Inspection of the proof of Proposition 2.4 shows that (23) continues to hold if the estimator $\hat{\tilde{\Phi}}_{n,P}(\cdot)$ is replaced by any of the estimators $\hat{\tilde{\Phi}}_{n,p}(\cdot)$ for $\mathcal{O}\leq p\leq P$. We also note that in case $\mathcal{O}=0$ the assumption of Proposition 2.4 is never satisfied (cf. Proposition 4.4 below), and hence, Theorem 2.3 always applies in that case. Furthermore, the case to which Proposition 2.4 applies is quite exceptional. In fact, under the assumptions of this proposition, the restricted estimators $A\tilde{\theta}(q)$ for $q\geq\mathcal{O}$ perform asymptotically as well as the unrestricted estimator $A\tilde{\theta}(P)$. This is again a consequence of Proposition 4.4.

We conclude this section by illustrating the above results with some important examples.

EXAMPLE 1 (The conditional distribution of $\tilde{\chi}$). Consider the model given in (1) with $\chi$ representing the parameter of interest. Using the general notation of Section 2, this corresponds to the case $A\theta = (\theta_1, \ldots, \theta_{\mathcal{O}})' = \chi$ with $A$ representing the $\mathcal{O} \times P$ matrix $(I_{\mathcal{O}} : 0)$. Here $k = \mathcal{O} > 0$. The c.d.f. $G_{n,\theta,\sigma}(\cdot|p)$ then represents the c.d.f. of $\sqrt{n}(\tilde{\chi} - \chi)$, conditional on the event $\hat{p} = p$. Assume first that $\lim_{n \to \infty} V'W/n \neq 0$. Then $C_\infty^{(r)} \neq 0$ holds for some $r > \mathcal{O}$. Consequently, the "impossibility" results for the estimation of $G_{n,\theta,\sigma}(t|\hat{p})$ given in Theorem 2.3 always apply. Next assume that $\lim_{n \to \infty} V'W/n = 0$. Then $C_\infty^{(r)} = 0$ for every $r > \mathcal{O}$. In this case Proposition 2.4 applies and a uniformly consistent estimator of $G_{n,\theta,\sigma}(t|\hat{p})$ indeed exists. Summarizing, we note that any estimator of $G_{n,\theta,\sigma}(t|\hat{p})$ suffers from the nonuniformity phenomenon, except in the special case where the columns of $V$ and $W$ are asymptotically orthogonal in the sense that $\lim_{n \to \infty} V'W/n = 0$. But this is precisely the situation where inclusion or exclusion of the regressors in $W$ has no effect on the (conditional) distribution of the estimator $\tilde{\chi}$ asymptotically; hence, it is not surprising that also the model selection procedure does not have an effect on the estimation of the c.d.f. of the post-model-selection estimator $\tilde{\chi}$. This observation may tempt one to enforce orthogonality between the columns of $V$ and $W$ by either replacing the columns of $V$ by their residuals from the projection on the column space of $W$ or vice versa. However, this is not helpful for the following reasons: In the first case one then in fact avoids model selection as all the restricted least-squares estimators for $\chi$ under consideration (and hence, also the post-model selection estimator $\tilde{\chi}$) in the reparameterized model coincide with the unrestricted least-squares estimator. In the second case the coefficients of the columns of $V$ in the reparameterized model no longer coincide with the parameter of interest $\chi$ (and again are estimated by one and the same estimator regardless of inclusion/exclusion of columns of the transformed $W$-matrix).

EXAMPLE 2 (The conditional distribution of $\tilde{\theta}$). For $A$ equal to $I_P$, the c.d.f. $G_{n,\theta,\sigma}(t|p)$ is the conditional c.d.f. of $\sqrt{n}(\tilde{\theta} - \theta)$ given $\hat{p} = p$. Here, $A\tilde{\theta}(q)$ reduces to $\tilde{\theta}(q)$, and hence, $A\tilde{\theta}(q)$ and $\tilde{\theta}_q(q)$ are perfectly correlated for every $q > \mathcal{O}$. Consequently, the "impossibility" result for estimation of $G_{n,\theta,\sigma}(t|\hat{p})$ given in Theorem 2.3 applies. We therefore see that estimation of the conditional distribution of the post-model-selection estimator of the entire parameter vector is always plagued by the nonuniformity phenomenon.

EXAMPLE 3 (The conditional distribution of a linear predictor). Suppose $A \neq 0$ is a $1 \times P$ vector and one is interested in estimating the conditional c.d.f. $G_{n,\theta,\sigma}(t|\hat{p})$ of the linear predictor $A\tilde{\theta}$. Then Theorem 2.3 and the discussion following Proposition 2.4 show that the nonuniformity phenomenon always arises in this estimation problem in case $\mathcal{O} = 0$. In case $\mathcal{O} > 0$, the nonuniformity problem is generically also present, except in the degenerate case where $C_\infty^{(q)} = 0$, for all $q$ satisfying $\mathcal{O} < q \leq P$ (in which case Proposition 4.4 shows that the least-squares predictors from all models $M_p$, $\mathcal{O} \leq p \leq P$, perform asymptotically equally well).

**3. Extensions to other model selection procedures including AIC.** In this section we show that the "impossibility" result obtained in the previous section for a "general-to-specific" model selection procedure carries over to a large class of model selection procedures, including Akaike's widely used AIC. Again, consider the linear regression model (5) with the same assumptions on the regressors and the errors as in Section 2. Let $\{0,1\}^P$ denote the set of all 0–1 sequences of length $P$. For each $\mathfrak{r} \in \{0,1\}^P$, let $M_\mathfrak{r}$ denote the set $\{\theta \in \mathbf{R}^P : \theta_i(1-\mathfrak{r}_i) = 0 \text{ for } 1 \leq i \leq P\}$, where $\mathfrak{r}_i$ represents the $i$th component of $\mathfrak{r}$. That is, $M_\mathfrak{r}$ describes a linear submodel with those parameters $\theta_i$ for which $\mathfrak{r}_i = 0$ restricted to zero. Now let $\mathfrak{R}$ be a user-supplied subset of $\{0,1\}^P$. We consider model selection procedures that select from the set $\mathfrak{R}$, or, equivalently, from the set of models $\{M_\mathfrak{r} : \mathfrak{r} \in \mathfrak{R}\}$. Note that there is now no assumption that the candidate models are nested (e.g., if $\mathfrak{R} = \{0,1\}^P$, all possible submodels are candidates for selection). Also, cases where the inclusion of a subset of regressors is undisputed on a priori grounds are obviously covered by this framework upon suitable choice of $\mathfrak{R}$.

We shall assume throughout this section that $\mathfrak{R}$ contains $\mathfrak{r}_{\text{full}} = (1,\ldots,1)$ and also at least one element $\mathfrak{r}_*$ satisfying $|\mathfrak{r}_*| = P - 1$, where $|\mathfrak{r}_*|$ represents the number of nonzero coordinates of $\mathfrak{r}_*$. Let $\hat{\mathfrak{r}}$ be an arbitrary model selection procedure, that is, $\hat{\mathfrak{r}} = \hat{\mathfrak{r}}(Y, X)$ is a random variable taking its values in $\mathfrak{R}$. We furthermore assume throughout this section that the model selection procedure $\hat{\mathfrak{r}}$ satisfies the following mild condition: For every $\mathfrak{r}_* \in \mathfrak{R}$ with $|\mathfrak{r}_*| = P - 1$, there exists a positive finite constant $c$ (possibly depending on $\mathfrak{r}_*$) such that, for every $\theta \in M_{\mathfrak{r}_*}$ which has exactly $P-1$ nonzero coordinates,

$$
\begin{aligned}
&\lim_{n \to \infty} P_{n,\theta,\sigma}(\{\hat{\mathfrak{r}} = \mathfrak{r}_{\text{full}}\} \blacktriangle \{|T_{\mathfrak{r}_*}| \geq c\}) \\
&\quad = \lim_{n \to \infty} P_{n,\theta,\sigma}(\{\hat{\mathfrak{r}} = \mathfrak{r}_*\} \blacktriangle \{|T_{\mathfrak{r}_*}| < c\}) = 0
\end{aligned}
\tag{24}
$$

holds for every $0 < \sigma < \infty$. Here $\blacktriangle$ denotes the symmetric difference operator and $T_{\mathfrak{r}_*}$ represents the usual $t$-statistic for testing the hypothesis $\theta_{i(\mathfrak{r}_*)} = 0$ in the full model, where $i(\mathfrak{r}_*)$ denotes the index of the unique coordinate of $\mathfrak{r}_*$ that equals zero.

The above condition is quite natural for the following reason: For $\theta \in M_{\mathfrak{r}_*}$ with exactly $P - 1$ nonzero coordinates, every reasonable model selection procedure will—with probability approaching unity—decide only between $M_{\mathfrak{r}_*}$ and $M_{\mathfrak{r}_{\text{full}}}$; it is then quite natural that this decision will be based (at least asymptotically) on the likelihood ratio between these two models, which in turn boils down to the $t$-statistic. As will be shown below, condition (24) holds in particular for AIC-like procedures.

Let $A$ be a nonstochastic $k \times P$ matrix of full row rank $k$, $1 \leq k \leq P$, as in Section 2.1. For every $\mathfrak{r} \in \mathfrak{R}$, we then consider the conditional c.d.f.

$$(25) \quad K_{n,\theta,\sigma}(t|\mathfrak{r}) = P_{n,\theta,\sigma}(\sqrt{n}A(\bar{\theta} - \theta) \leq t|\hat{\mathfrak{r}} = \mathfrak{r}) \qquad (t \in \mathbf{R}^k)$$

of a linear transformation of the post-model-selection estimator $\bar{\theta}$ obtained from the model selection procedure $\hat{\mathfrak{r}}$, that is,

$$\bar{\theta} = \sum_{\mathfrak{r} \in \mathfrak{R}} \tilde{\theta}(\mathfrak{r})\mathbf{1}(\hat{\mathfrak{r}} = \mathfrak{r}),$$

where the $P \times 1$ vector $\tilde{\theta}(\mathfrak{r})$ represents the restricted least-squares estimator obtained from model $M_\mathfrak{r}$, with the convention that $\tilde{\theta}(\mathfrak{r}) = 0 \in \mathbf{R}^P$ in case $\mathfrak{r} = (0, \ldots, 0)$. [In case $P_{n,\theta,\sigma}(\hat{\mathfrak{r}} = \mathfrak{r}) = 0$, we define $K_{n,\theta,\sigma}(t|p)$ equal to, say, the c.d.f. of point-mass at zero in $\mathbf{R}^k$. This is done just for the sake of definiteness and has no effect on the results given below. For most model selection procedures, the probability $P_{n,\theta,\sigma}(\hat{\mathfrak{r}} = \mathfrak{r})$ will be positive for any $\mathfrak{r} \in \mathfrak{R}$ anyway.] We also introduce

$$(26) \quad K_{n,\theta,\sigma}(t|\hat{\mathfrak{r}}) = \sum_{\mathfrak{r} \in \mathfrak{R}} K_{n,\theta,\sigma}(t|\mathfrak{r})\mathbf{1}(\hat{\mathfrak{r}} = \mathfrak{r}) \qquad (t \in \mathbf{R}^k).$$

We then obtain the following result for estimation of $K_{n,\theta,\sigma}(t|\hat{\mathfrak{r}})$ at a *fixed* value of the argument $t$ which parallels the corresponding "impossibility" result in Section 2.2.2.

THEOREM 3.1. *Let $\mathfrak{r}_* \in \mathfrak{R}$ satisfy $|\mathfrak{r}_*| = P - 1$, and let $i(\mathfrak{r}_*)$ denote the index of the unique coordinate of $\mathfrak{r}_*$ that equals zero; furthermore, let $c$ be the constant in (24) corresponding to $\mathfrak{r}_*$. Suppose that $A\tilde{\theta}(\mathfrak{r}_{\text{full}})$ and $\tilde{\theta}_{i(\mathfrak{r}_*)}(\mathfrak{r}_{\text{full}})$ are asymptotically correlated, that is, $AQ^{-1}e_{i(\mathfrak{r}_*)} \neq 0$, where $e_{i(\mathfrak{r}_*)}$ denotes the $i(\mathfrak{r}_*)$th standard basis vector in $\mathbf{R}^P$. Then for every $\theta \in M_{\mathfrak{r}_*}$ which has exactly $P - 1$ nonzero coordinates, for each $0 < \sigma < \infty$ and for each $t \in \mathbf{R}^k$, the following holds: There exist $\delta_0 > 0$ and $0 < \rho_0 < \infty$ such that any estimator $\hat{K}_n(t|\hat{\mathfrak{r}})$ of $K_{n,\theta,\sigma}(t|\hat{\mathfrak{r}})$ satisfying*

$$(27) \quad P_{n,\theta,\sigma}(|\hat{K}_n(t|\hat{\mathfrak{r}}) - K_{n,\theta,\sigma}(t|\hat{\mathfrak{r}})| > \delta) \overset{n\to\infty}{\longrightarrow} 0$$

*for each $\delta > 0$ (in particular, every estimator that is consistent) also satisfies*

$$(28) \quad \liminf_{n\to\infty} \sup_{\substack{\vartheta \in \mathbf{R}^P \\ \|\vartheta - \theta\| < \rho_0/\sqrt{n}}} P_{n,\vartheta,\sigma}(|\hat{K}_n(t|\hat{\mathfrak{r}}) - K_{n,\vartheta,\sigma}(t|\hat{\mathfrak{r}})| > \delta_0)$$
$$\geq 2(1 - \Phi(c)) > 0.$$

*The constants $\delta_0$ and $\rho_0$ may be chosen in such a way that they depend only on $t, Q, A, \sigma$ and $c$. Moreover,*

$$\text{(29)} \quad \liminf_{n \to \infty} \inf_{\hat{K}_n(t|\hat{\mathfrak{r}})} \sup_{\substack{\vartheta \in \mathbf{R}^P \\ \|\vartheta - \theta\| < \rho_0/\sqrt{n}}} P_{n,\vartheta,\sigma}(|\hat{K}_n(t|\hat{\mathfrak{r}}) - K_{n,\vartheta,\sigma}(t|\hat{\mathfrak{r}})| > \delta_0) > 0$$

*and*

$$\text{(30)} \quad \sup_{\delta > 0} \liminf_{n \to \infty} \inf_{\hat{K}_n(t|\hat{\mathfrak{r}})} \sup_{\substack{\vartheta \in \mathbf{R}^P \\ \|\vartheta - \theta\| < \rho_0/\sqrt{n}}} P_{n,\vartheta,\sigma}(|\hat{K}_n(t|\hat{\mathfrak{r}}) - K_{n,\vartheta,\sigma}(t|\hat{\mathfrak{r}})| > \delta)$$

$$\geq 1 - \Phi(c) > 0$$

*hold, where the infima in* (29) *and* (30) *extend over all estimators $\hat{K}_n(t|\hat{\mathfrak{r}})$ of $K_{n,\theta,\sigma}(t|\hat{\mathfrak{r}})$. [The lower bound in* (28) *is nothing other than $\lim_{n \to \infty} P_{n,\theta,\sigma}(\hat{\mathfrak{r}} = \mathfrak{r}_{\mathrm{full}})$.]*

The basic condition (24) on the model selection procedure employed in the above results will certainly hold for any hypothesis testing procedure that (i) asymptotically selects only correct models, (ii) employs a likelihood ratio test (or an asymptotically equivalent test) for testing $M_{\mathfrak{r}_{\mathrm{full}}}$ versus smaller models [at least versus the models $M_{\mathfrak{r}_*}$ with $\mathfrak{r}_*$ as in condition (24)], and (iii) uses a critical value for the likelihood ratio test that converges to a finite positive constant. In particular, this applies to usual thresholding procedures, as well as to a variant of the "general-to-specific" procedure discussed in Section 2 where the error variance in the construction of the test statistic for hypothesis $H_0^p$ is estimated from the fitted model $M_p$ rather than from the overall model. We next verify condition (24) for AIC-like procedures. Let $RSS(\mathfrak{r})$ denote the residual sum of squares from the regression employing model $M_\mathfrak{r}$ and set

$$\text{(31)} \quad IC(\mathfrak{r}) = \log(RSS(\mathfrak{r})) + |\mathfrak{r}|\Upsilon_n/n,$$

where $\Upsilon_n \geq 0$ denotes a sequence of real numbers satisfying $\lim_{n \to \infty} \Upsilon_n = \Upsilon$ and $\Upsilon$ is a positive real number. Of course, $IC(\mathfrak{r}) = AIC(\mathfrak{r})$ if $\Upsilon_n = 2$. The model selection procedure $\hat{\mathfrak{r}}_{IC}$ is then defined as a minimizer (more precisely, as a measurable selection from the set of minimizers) of $IC(\mathfrak{r})$ over $\mathfrak{R}$. It is well known that the probability that $\hat{\mathfrak{r}}_{IC}$ selects an incorrect model converges to zero. Hence, elementary calculations show that condition (24) is satisfied for $c = \Upsilon^{1/2}$.

## 4. Further theoretical results.

4.1. *"General-to-specific" model selection procedure.* In this section we provide "impossibility" results for estimation of $G_{n,\theta,\sigma}(t|p)$ for *given* $p$ which are parallel to the "impossibility" result for estimation of $G_{n,\theta,\sigma}(t|\hat{p})$ given in Section 2.2.2. The results presented below can be viewed as conditional counterparts to the results in that earlier section. Apart from being of interest on their own, the results given below also form the essential building blocks for the "impossibility" result in Section 2.2.2. In the next two theorems we shall consider estimation of $G_{n,\theta,\sigma}(t|p)$ at a *fixed* value of the argument $t$.

THEOREM 4.1. *Let $p$ satisfy $\mathcal{O} < p \leq P$. Suppose that $A\tilde{\theta}(p)$ and $\tilde{\theta}_p(p)$ are asymptotically correlated, that is, $C_\infty^{(p)} \neq 0$. Then the following holds for each $\theta \in M_{p-1}$, $0 < \sigma < \infty$, and for each $t \in \mathbf{R}^k$ with the property that the set $\{z \in \mathbf{R}^p : A[p]z \leq t\}$ has positive Lebesgue measure in $\mathbf{R}^p$:*

(a) *There exist $\delta_0 > 0$ and $0 < \rho_0 < \infty$ such that any estimator $\hat{G}_n(t|p)$ for $G_{n,\theta,\sigma}(t|p)$ satisfying*

$$P_{n,\theta,\sigma}(|\hat{G}_n(t|p) - G_{n,\theta,\sigma}(t|p)| > \delta) \stackrel{n \to \infty}{\longrightarrow} 0 \tag{32}$$

*for each $\delta > 0$ (in particular, every estimator that is consistent over $M_p$) also satisfies*

$$\sup_{\substack{\vartheta \in M_p \\ \|\vartheta - \theta\| < \rho_0/\sqrt{n}}} P_{n,\vartheta,\sigma}(|\hat{G}_n(t|p) - G_{n,\vartheta,\sigma}(t|p)| > \delta_0) \stackrel{n \to \infty}{\longrightarrow} 1. \tag{33}$$

*The constants $\delta_0$ and $\rho_0$ may be chosen in such a way that they depend only on $t, Q, A, \sigma$ and the critical value $c_p$. Moreover,*

$$\liminf_{n \to \infty} \inf_{\hat{G}_n(t|p)} \sup_{\substack{\vartheta \in M_p \\ \|\vartheta - \theta\| < \rho_0/\sqrt{n}}} P_{n,\vartheta,\sigma}(|\hat{G}_n(t|p) - G_{n,\vartheta,\sigma}(t|p)| > \delta_0) > 0 \tag{34}$$

*and*

$$\sup_{\delta > 0} \liminf_{n \to \infty} \inf_{\hat{G}_n(t|p)} \sup_{\substack{\vartheta \in M_p \\ \|\vartheta - \theta\| < \rho_0/\sqrt{n}}} P_{n,\vartheta,\sigma}(|\hat{G}_n(t|p) - G_{n,\vartheta,\sigma}(t|p)| > \delta) \geq \tfrac{1}{2} \tag{35}$$

*hold, where the infima in* (34) *and* (35) *extend over all estimators $\hat{G}_n(t|p)$ of $G_{n,\theta,\sigma}(t|p)$.*

(b) *The above continues to hold with $P_{n,\cdot,\sigma}(\cdot|\hat{p} = p)$ replacing $P_{n,\cdot,\sigma}(\cdot)$.*

The condition on the set $\{z \in \mathbf{R}^p : A[p]z \leq t\}$ in the above theorem is easily seen to be equivalent to the condition that $A[p]z < t$ holds for some $z \in \mathbf{R}^p$. It is trivially always satisfied whenever $t > 0$. The condition on $\{z \in \mathbf{R}^p : A[p]z \leq t\}$ is certainly satisfied for *every* $t \in \mathbf{R}^k$ if the matrix $A[p]$ has full row rank $k$. We shall repeatedly use the observation that the latter rank condition is always met if $p > 0$ is the maximal order for which $C_\infty^{(p)} \neq 0$ holds. [This follows from Proposition 4.4(a), (c).]

As a point of interest, we note that the nonuniformity phenomenon described in Theorem 4.1 occurs *within* the model $M_p$, which contains only parameters for which the selected model is correct; that is, in (33)–(35) the suprema with respect to $\vartheta$ extend only over subsets of $M_p$. That is, it is typically even impossible to construct an estimator of $G_{n,\theta,\sigma}(t|p)$ which performs satisfactorily for those local perturbations $\vartheta$ of the true parameter $\theta \in M_{p-1}$ for which the selected model is correct.

Consider next the case where Theorem 4.1 does not apply, that is, the model order $p$ under consideration is such that either $p = \mathcal{O}$, or $p > \mathcal{O}$ but $C_\infty^{(p)} = 0$, or $p > \mathcal{O}$ and $C_\infty^{(p)} \neq 0$ but the set $\{z \in \mathbf{R}^p : A[p]z \leq t\}$ has Lebesgue measure zero. In that case, it is indeed possible to construct an estimator of $G_{n,\theta,\sigma}(t|p)$ that is uniformly consistent over $\theta \in M_p$. However, this result provides little consolation, because the uniform consistency over $\theta \in M_p$ typically breaks down already in $1/\sqrt{n}$-"neighborhoods" of $M_p$, and results analogous to (33)–(35) can be established over such neighborhoods, even if Theorem 4.1 does not apply. This is of relevance as true parameter values in such $1/\sqrt{n}$-"neighborhoods" result in a positive probability of selecting the model $M_p$; see Proposition A.2 in Appendix A.

THEOREM 4.2. *Let $p$ satisfy $\mathcal{O} \leq p < P$. Suppose that $A\tilde{\theta}(q)$ and $\tilde{\theta}_q(q)$ are asymptotically correlated, that is, $C_\infty^{(q)} \neq 0$, for some $q$ satisfying $p < q \leq P$, and let $q^*$ denote the largest $q$ with this property. Then the following holds for each $\theta \in M_p$, $0 < \sigma < \infty$, and for each $t \in \mathbf{R}^k$:*

(a) *There exist $\delta_0 > 0$ and $0 < \rho_0 < \infty$ such that any estimator $\hat{G}_n(t|p)$ for $G_{n,\theta,\sigma}(t|p)$ satisfying*

$$(36) \qquad P_{n,\theta,\sigma}(|\hat{G}_n(t|p) - G_{n,\theta,\sigma}(t|p)| > \delta) \overset{n \to \infty}{\longrightarrow} 0$$

*for each $\delta > 0$ (in particular, every estimator that is consistent over $M_p$) also satisfies*

$$(37) \qquad \sup_{\substack{\vartheta \in M_{q^*} \\ \|\vartheta - \theta\| < \rho_0/\sqrt{n}}} P_{n,\vartheta,\sigma}(|\hat{G}_n(t|p) - G_{n,\vartheta,\sigma}(t|p)| > \delta_0) \overset{n \to \infty}{\longrightarrow} 1.$$

*The constants $\delta_0$ and $\rho_0$ may be chosen in such a way that they depend only on $t, Q, A, \sigma$ and the critical value $c_p$. Moreover,*

$$(38) \qquad \liminf_{n \to \infty} \inf_{\hat{G}_n(t|p)} \sup_{\substack{\vartheta \in M_{q^*} \\ \|\vartheta - \theta\| < \rho_0/\sqrt{n}}} P_{n,\vartheta,\sigma}(|\hat{G}_n(t|p) - G_{n,\vartheta,\sigma}(t|p)| > \delta_0) > 0$$

*and*

(39) $$\sup_{\delta>0} \liminf_{n\to\infty} \inf_{\hat{G}_n(t|p)} \sup_{\substack{\vartheta\in M_{q^*} \\ \|\vartheta-\theta\|<\rho_0/\sqrt{n}}} P_{n,\vartheta,\sigma}(|\hat{G}_n(t|p) - G_{n,\vartheta,\sigma}(t|p)| > \delta) \geq \tfrac{1}{2}$$

*hold, where the infima in* (38) *and* (39) *extend over* all *estimators* $\hat{G}_n(t|p)$ *of* $G_{n,\theta,\sigma}(t|p)$.

(b) *The above continues to hold with* $P_{n,\cdot,\sigma}(\cdot|\hat{p}=p)$ *replacing* $P_{n,\cdot,\sigma}(\cdot)$.

We stress here once more that the probability of selecting the model order $p$ is bounded away from zero uniformly over the $\vartheta$-sets appearing in the suprema in (37)–(39); see Proposition A.2 in Appendix A. Hence, the nonuniformity phenomenon we observe is not an artifact resulting from conditioning on an unlikely event. It is also worth noting that the lower bounds in the above results are as large as 1 and 1/2, respectively.

Summarizing so far, we see that it is impossible to construct an estimator of $G_{n,\theta,\sigma}(t|p)$ which performs reasonably well in a neighborhood of the true parameter $\theta$ ($\theta \in M_p$), whenever the model order $p$ considered has the property that $A\tilde{\theta}(q)$ and $\tilde{\theta}_q(q)$ are asymptotically correlated for some $q$ with $\max\{p, \mathcal{O}+1\} \leq q \leq P$, as then either Theorem 4.1 or Theorem 4.2 applies. In particular, no uniformly consistent estimator exists, not even locally. In the remaining case, that is, when $A\tilde{\theta}(q)$ and $\tilde{\theta}_q(q)$ are asymptotically uncorrelated for each $q$ in the range $\max\{p, \mathcal{O}+1\} \leq q \leq P$, it is indeed possible to construct an estimator of $G_{n,\theta,\sigma}(t|p)$ which is uniformly consistent (even in the total variation distance) over $1/\sqrt{n}$-"neighborhoods" of $M_p$, as shown next.

PROPOSITION 4.3. *Let $p$ satisfy $\mathcal{O} \leq p \leq P$. Suppose that $A\tilde{\theta}(q)$ and $\tilde{\theta}_q(q)$ are asymptotically uncorrelated, that is, $C_\infty^{(q)} = 0$ for each $q = \max\{p, \mathcal{O}+1\}, \ldots, P$. Then*

(40) $$\sup_{\substack{\theta\in\mathbf{R}^P \\ \|\theta[\neg p]\|<\rho/\sqrt{n}}} \sup_{\substack{\sigma\in\mathbf{R} \\ \sigma_*\leq\sigma\leq\sigma^*}} P_{n,\theta,\sigma}(\|\hat{\tilde{\Phi}}_{n,p}(\cdot) - G_{n,\theta,\sigma}(\cdot|p)\|_{\mathrm{TV}} > \delta) \overset{n\to\infty}{\longrightarrow} 0$$

*holds for each $\delta > 0$, for each $0 < \rho < \infty$, and for any constants $\sigma_*$ and $\sigma^*$ satisfying $0 < \sigma_* \leq \sigma^* < \infty$. The result* (40) *also holds with $P_{n,\theta,\sigma}(\cdot|\hat{p}=p)$ replacing $P_{n,\theta,\sigma}(\cdot)$. [In case $p=P$, the first supremum in* (40) *is to be interpreted as extending over all $\theta \in \mathbf{R}^P$. Furthermore, the case $p=0$ is impossible in view of Proposition* 4.4 *below.*]

If the uncorrelatedness assumptions in the proposition even hold for all finite $n$, then the c.d.f. $G_{n,\theta,\sigma}(\cdot|p)$ can be seen to reduce to the normal c.d.f.

$\Phi_{n,p}(\cdot)$ and, hence, can be estimated uniformly consistently over the larger space $M_P \times [\sigma_*, \sigma^*]$.

Clearly, the case to which Proposition 4.3 applies is quite exceptional. In fact, under the assumptions of this proposition, the restricted estimators $A\tilde{\theta}(q)$ for $q \geq \max\{p-1, \mathcal{O}\}$ perform asymptotically as well as the unrestricted estimator $A\tilde{\theta}(P)$. This is a consequence of the following result.

PROPOSITION 4.4. *Let $p$ satisfy $0 < p < P$. Then the following statements are equivalent:*

(a) $A\tilde{\theta}(q)$ *and* $\tilde{\theta}_q(q)$ *are asymptotically uncorrelated, that is,* $C_\infty^{(q)} = (0, \ldots, 0)'$ *for each* $q = p+1, \ldots, P$.

(b) $A\tilde{\theta}(p)$ *is an asymptotically unbiased estimator of* $A\theta$ *($\theta \in \mathbf{R}^P$)*.

(c) *The asymptotic variance–covariance matrices of* $\sqrt{n}A\tilde{\theta}(p)$ *and* $\sqrt{n}A\tilde{\theta}(P)$ *are identical.*

*In case $p = P$, the above statements are always trivially satisfied. In case $p = 0$, these statements are never satisfied.*

It is easy to see that any of the above statements is equivalent to asymptotic unbiasedness of $A\tilde{\theta}(q)$ for all $q = p, \ldots, P$, and further, also is equivalent to all the asymptotic variance–covariance matrices of $\sqrt{n}A\tilde{\theta}(q)$ for $q = p, \ldots, P$ being identical. Furthermore, a finite sample version of Proposition 4.4 can also easily be derived from the discussion following (19) in [10]. In fact, it is shown in that reference for any given sample size that uncorrelatedness of $A\tilde{\theta}(q)$ and $\tilde{\theta}_q(q)$ for $q = p+1, \ldots, P$ is equivalent to the estimators $A\tilde{\theta}(p)$ and $A\tilde{\theta}(P)$ being identical, which, in turn, is equivalent to all the estimators $A\tilde{\theta}(q)$ being identical for $q = p, \ldots, P$.

We conclude this section by illustrating the above results with some important examples.

EXAMPLE 1 CONTINUED (The conditional distribution of $\tilde{\chi}$). Assume first that $\lim_{n \to \infty} V'W/n \neq 0$ is satisfied. Then, as already noted, $C_\infty^{(r)} \neq 0$ holds for some $r > \mathcal{O}$. Consequently, for any such $r$, the "impossibility" results in Theorem 4.1 apply with $p = r$ (observe that $\mathrm{rank}(A[p]) = \mathcal{O} = k$ always holds for $p = r > \mathcal{O}$ and, hence, the condition on $t$ in that theorem is always satisfied). Furthermore, the "impossibility" results in Theorem 4.2 apply for any $p$ satisfying $\mathcal{O} \leq p < r$ for some $r$ as above. Next assume that $\lim_{n \to \infty} V'W/n = 0$. Then $C_\infty^{(r)} = 0$ for every $r > \mathcal{O}$. In this case Proposition 4.3 applies for every $\mathcal{O} \leq p \leq P$, and an estimator of $G_{n,\theta,\sigma}(t|p)$ that is uniformly consistent over $1/\sqrt{n}$-"neighborhoods" of $M_p$ indeed exists.

EXAMPLE 2 CONTINUED (The conditional distribution of $\tilde{\theta}$). As already noted, we have here $A = I_P$ and $A\tilde{\theta}(q)$ and $\tilde{\theta}_q(q)$ are perfectly correlated for every $q > \mathcal{O}$. Therefore, Theorem 4.1 applies for all $t \in \mathbf{R}^k$ if $p = P$, and Theorem 4.2 applies in case $p < P$. (In the latter case, Theorem 4.1 still applies for certain $t \in \mathbf{R}^k$.) Consequently, estimation of the conditional distribution $G_{n,\theta,\sigma}(t|p)$ of the entire parameter vector is always plagued by the nonuniformity phenomenon.

EXAMPLE 4 (The conditional distribution of the unrestricted components of $\tilde{\theta}$). Let $r > 0$ be a given model order. Conditional on the event $\hat{p} = r$, the last $P - r$ components of $\tilde{\theta}$ are restricted to zero. If $A$ is the $r \times P$ matrix $(I_r : 0)$, then the c.d.f. $G_{n,\theta,\sigma}(t|r)$ is the conditional c.d.f. of the first $r$ (unrestricted) components of $\sqrt{n}(\tilde{\theta} - \theta)$ given the event $\hat{p} = r$. In this case $A\tilde{\theta}(r)$ and $\tilde{\theta}_r(r)$ are perfectly correlated. If $r > \mathcal{O}$, Theorem 4.1 immediately applies with $p = r$, because $\text{rank}(A[r]) = r$ entails that the condition on $t$ in that theorem is always satisfied. In case $\mathcal{O} \leq r < P$ and $\lim_{n\to\infty} X[r]'X[\neg r]/n \neq 0$, Theorem 4.2 applies with $p = r$, since, under the latter condition on the regressors, $C_\infty^{(q)} \neq 0$ holds for some $q > r$. As a consequence, the nonuniformity phenomenon is always present when estimating this conditional c.d.f., except in the very special case where $r = \mathcal{O} > 0$ and $\lim_{n\to\infty} X[r]'X[\neg r]/n = 0$ simultaneously hold; in this case Proposition 4.3 applies.

4.2. *Other model selection procedures including AIC.* We use the notation and assumptions of Section 3 here. In particular, the model selection procedure $\hat{\mathfrak{r}}$ is assumed to satisfy condition (24). The proof of Theorem 3.1 relies on the subsequent result, which is of interest also in itself. Similarly as in the preceding sections, estimation of $K_{n,\theta,\sigma}(t|\mathfrak{r})$ at a *fixed* value of the argument $t$ is considered.

THEOREM 4.5. *Let $\mathfrak{r}_* \in \mathfrak{R}$ satisfy $|\mathfrak{r}_*| = P - 1$, and let $i(\mathfrak{r}_*)$ denote the index of the unique coordinate of $\mathfrak{r}_*$ that equals zero; furthermore, let c be the constant in* (24) *corresponding to $\mathfrak{r}_*$. Suppose that $A\tilde{\theta}(\mathfrak{r}_{\text{full}})$ and $\tilde{\theta}_{i(\mathfrak{r}_*)}(\mathfrak{r}_{\text{full}})$ are asymptotically correlated, that is, $AQ^{-1}e_{i(\mathfrak{r}_*)} \neq 0$, where $e_{i(\mathfrak{r}_*)}$ denotes the $i(\mathfrak{r}_*)$th standard basis vector in $\mathbf{R}^P$. Then for every $\theta \in M_{\mathfrak{r}_*}$ which has exactly $P - 1$ nonzero coordinates, for each $0 < \sigma < \infty$, and for each $t \in \mathbf{R}^k$, the following holds with $\mathfrak{r} = \mathfrak{r}_*$ as well as with $\mathfrak{r} = \mathfrak{r}_{\text{full}}$:*

(a) *There exist $\delta_0 > 0$ and $0 < \rho_0 < \infty$ such that any estimator $\hat{K}_n(t|\mathfrak{r})$ for $K_{n,\theta,\sigma}(t|\mathfrak{r})$ satisfying*

$$(41) \qquad P_{n,\theta,\sigma}(|\hat{K}_n(t|\mathfrak{r}) - K_{n,\theta,\sigma}(t|\mathfrak{r})| > \delta) \stackrel{n\to\infty}{\longrightarrow} 0$$

*for each $\delta > 0$ (in particular, every estimator that is consistent over $M_{\mathfrak{r}}$) also satisfies*

$$(42) \quad \sup_{\substack{\vartheta \in \mathbf{R}^P \\ \|\vartheta - \theta\| < \rho_0/\sqrt{n}}} P_{n,\vartheta,\sigma}(|\hat{K}_n(t|\mathfrak{r}) - K_{n,\vartheta,\sigma}(t|\mathfrak{r})| > \delta_0) \overset{n \to \infty}{\longrightarrow} 1.$$

*The constants $\delta_0$ and $\rho_0$ may be chosen in such a way that they depend only on $t, Q, A, \sigma$ and also on $c$ in case $\mathfrak{r} = \mathfrak{r}_{\text{full}}$. Moreover,*

$$(43) \quad \liminf_{n \to \infty} \inf_{\hat{K}_n(t|\mathfrak{r})} \sup_{\substack{\vartheta \in \mathbf{R}^P \\ \|\vartheta - \theta\| < \rho_0/\sqrt{n}}} P_{n,\vartheta,\sigma}(|\hat{K}_n(t|\mathfrak{r}) - K_{n,\vartheta,\sigma}(t|\mathfrak{r})| > \delta_0) > 0$$

*and*

$$(44) \quad \sup_{\delta > 0} \liminf_{n \to \infty} \inf_{\hat{K}_n(t|\mathfrak{r})} \sup_{\substack{\vartheta \in \mathbf{R}^P \\ \|\vartheta - \theta\| < \rho_0/\sqrt{n}}} P_{n,\vartheta,\sigma}(|\hat{K}_n(t|\mathfrak{r}) - K_{n,\vartheta,\sigma}(t|\mathfrak{r})| > \delta) \geq \tfrac{1}{2}$$

*hold, where the infima in (43) and (44) extend over* all *estimators $\hat{K}_n(t|\mathfrak{r})$ of $K_{n,\theta,\sigma}(t|\mathfrak{r})$.*

  (b) *The above continues to hold with $P_{n,\cdot,\sigma}(\cdot|\hat{\mathfrak{r}} = \mathfrak{r})$ replacing $P_{n,\cdot,\sigma}(\cdot)$.*

We note that the conditional probability in Theorem 4.5(b) is eventually well defined; see (61)–(62) in Appendix E.

### 4.3. Remarks and extensions.

REMARK 4.6. Although not emphasized in the notation, all results in the paper also hold if the elements of the design matrix $X$ depend on sample size. Furthermore, all results are expressed solely in terms of the distributions $P_{n,\theta,\sigma}(\cdot)$ of $Y$, and hence, they also apply if the elements of $Y$ depend on sample size, including the case where the random vectors $Y$ are defined on different probability spaces for different sample sizes.

REMARK 4.7. The model selection procedure introduced in Section 2 is based on a sequence of tests which use critical values $c_p$ that do not depend on sample size and satisfy $0 < c_p < \infty$ for $\mathcal{O} < p \leq P$. If these critical values are allowed to depend on sample size such that they now satisfy $c_{n,p} \to c_{\infty,p}$ as $n \to \infty$ with $0 < c_{\infty,p} < \infty$ for $\mathcal{O} < p \leq P$, the results in [12], as well as in [10, 11], continue to hold; see Remark 6.2(i) in [12] and Remark 6.1(ii) in [10]. As a consequence, the results in the present paper can also be extended to this case quite easily.

REMARK 4.8. The "impossibility" results given in Theorems 2.3, 3.1, 4.1, 4.2 and 4.5 (as well as the variants thereof discussed in the subsequent

Remark [4.9](#)) also hold for the class of all randomized estimators (with $P^*_{n,\theta,\sigma}$ replacing $P_{n,\theta,\sigma}$ in those results, where $P^*_{n,\theta,\sigma}$ denotes the distribution of the randomized sample). This follows immediately from Lemma 3.6 and the attending discussion in [16].

REMARK 4.9. Results similar to the ones in Sections [2.2.2](#) and [4.1](#) can also be obtained for estimation of the asymptotic c.d.f. $G_{\infty,\theta,\sigma}(t|p)$. Since these results are of limited interest, we omit them. In particular, note that an "impossibility" result for estimation of $G_{\infty,\theta,\sigma}(t|p)$ per se does *not* imply a corresponding "impossibility" result for estimation of $G_{n,\theta,\sigma}(t|p)$, since $G_{n,\theta,\sigma}(t|p)$ does in general not converge uniformly to $G_{\infty,\theta,\sigma}(t|p)$ over the relevant subsets in the parameter space; see Remark 4.11 in [14]. (Appropriate analogues apply to the model selection procedures considered in Sections [3](#) and [4.2](#).)

**5. Conclusion.** Despite the fact that we have shown that consistent estimators for the conditional distribution of a post-model-selection estimator can be constructed with relative ease, we have also demonstrated that no estimator of this conditional distribution can have satisfactory performance (locally) uniformly in the parameter space, even asymptotically. In particular, no (locally) uniformly consistent estimator of the conditional distribution exists. Hence, the answer to the question posed in the title has to be negative. The results in the present paper also cover the case of linear functions (e.g., predictors) of the post-model-selection estimator. Corresponding results for the unconditional distribution of the post-model-selection estimator are presented in a companion paper [13].

The "impossibility" results are derived in the framework of a normal linear regression model (and a fortiori these results continue to hold in any model which includes the normal linear regression model as a special case). Furthermore, there is no reason to believe that the situation will get any better in more complex statistical models that allow, for example, for nonlinearity or dependent data. In fact, similar results can be obtained in general statistical models, for example, as long as standard regularity conditions for maximum likelihood theory are satisfied.

The results in the present paper are derived for a large class of conservative model selection procedures (i.e., procedures that select overparameterized models with positive probability asymptotically), including Akaike's AIC and typical "general-to-specific" hypothesis testing procedures. For consistent model selection procedures—such as BIC or testing procedures with suitably diverging critical values $c_p$ (cf. [2])—the (pointwise) asymptotic distribution is always normal. (This is elementary; cf. Lemma 1 in [18].) However, as discussed at length in [15], this asymptotic normality result paints a misleading picture of the finite sample distribution, which can be far from normal, the convergence of the finite-sample distribution to the asymptotic normal distribution not being uniform. "Impossibility" results similar to the ones presented here can also be obtained for post-model-selection estimators based on consistent model selection procedures. These will be discussed in detail elsewhere. For a simple special case, such an "impossibility" result is given in Section 2.3 of [16].

The "impossibility" of estimating the distribution of the post-model-selection estimator does not per se preclude the possibility of conducting valid inference after model selection, a topic that deserves further study. However, it certainly makes this a more challenging task.

## APPENDIX A: THE LARGE-SAMPLE LIMIT OF $G_{N,\theta,\sigma}(T|P)$

For $p$ satisfying $0 < p \leq P$, partition the matrix $Q = \lim_{n\to\infty} X'X/n$ as

$$Q = \begin{pmatrix} Q[p:p] & Q[p:\neg p] \\ Q[\neg p:p] & Q[\neg p:\neg p] \end{pmatrix},$$

where $Q[p:p]$ is a $p \times p$ matrix. For $p = 1, \ldots, P$, define

$$\begin{aligned}
\xi_{\infty,p}^2 &= (Q[p:p]^{-1})_{p,p}, \\
(45) \quad \zeta_{\infty,p}^2 &= \xi_{\infty,p}^2 - C_\infty^{(p)\prime}(A[p]Q[p:p]^{-1}A[p]')^- C_\infty^{(p)}, \\
b_{\infty,p} &= C_\infty^{(p)\prime}(A[p]Q[p:p]^{-1}A[p]')^-,
\end{aligned}$$

where $C_\infty^{(p)} = A[p]Q[p:p]^{-1}e_p$, with $e_p$ denoting the $p$th standard basis vector in $\mathbf{R}^p$; furthermore, take $\zeta_{\infty,p}$ and $\xi_{\infty,p}$ as the nonnegative square roots of $\zeta_{\infty,p}^2$ and $\xi_{\infty,p}^2$, respectively. As the notation suggests, $\Phi_{\infty,p}(t)$ is the large-sample limit of $\Phi_{n,p}(t)$, both defined in Section 2. Moreover, $C_\infty^{(p)}$, $\xi_{\infty,p}^2$ and $\zeta_{\infty,p}^2$ are the limits of $C_n^{(p)}$, $\xi_{n,p}^2$ and $\zeta_{n,p}^2$, respectively, and $b_{n,p}z$ converges to $b_{\infty,p}z$ for each $z$ in the column-space of $A[p]$. See Lemma A.2 in [10].

The next result is taken from Corollary 5.4 in [10] and describes the large-sample limit of the conditional c.d.f. under local alternatives to $\theta$, under the assumption that the selected model $M_p$ is a correct model for $\theta$. Recall that the total variation distance between two c.d.f.s $G$ and $G^*$ on $\mathbf{R}^k$ is defined as $\|G - G^*\|_{\mathrm{TV}} = \sup_E |G(E) - G^*(E)|$, where the supremum is taken over all Borel sets $E$. Clearly, the relation $|G(t) - G^*(t)| \leq \|G - G^*\|_{\mathrm{TV}}$ holds for all $t \in \mathbf{R}^k$. Thus, if $G$ and $G^*$ are close with respect to the total variation distance, then $G(t)$ is close to $G^*(t)$, uniformly in $t$.

PROPOSITION A.1. *Let $p$ satisfy $\mathcal{O} \leq p \leq P$. Suppose $\theta \in \mathbf{R}^P$ satisfies $\theta \in M_p$, that is, $p_0(\theta) \leq p$ holds. Moreover, let $\gamma \in \mathbf{R}^P$ and let $\sigma^{(n)}$ be a sequence of positive real numbers which converges to a (finite) limit $\sigma > 0$*

POST-MODEL-SELECTION ESTIMATORS 27

as $n \to \infty$. Then the conditional c.d.f. $G_{n,\theta+\gamma/\sqrt{n},\sigma^{(n)}}(t|p)$ converges to a limit $G_{\infty,\theta,\sigma,\gamma}(t|p)$ in total variation, that is,

$$(46) \qquad \|G_{n,\theta+\gamma/\sqrt{n},\sigma^{(n)}}(\cdot|p) - G_{\infty,\theta,\sigma,\gamma}(\cdot|p)\|_{\mathrm{TV}} \stackrel{n\to\infty}{\longrightarrow} 0.$$

The large-sample limit c.d.f. $G_{\infty,\theta,\sigma,\gamma}(t|p)$ is given as follows: In case $p = \max\{p_0(\theta), \mathcal{O}\}$,

$$(47) \qquad G_{\infty,\theta,\sigma,\gamma}(t|p) = \Phi_{\infty,p}(t - \beta^{(p)}).$$

Here,

$$\beta^{(p)} = A \begin{pmatrix} Q[p\!:\!p]^{-1} Q[p\!:\!\neg p]\gamma[\neg p] \\ -\gamma[\neg p] \end{pmatrix} \qquad (0 \le p \le P),$$

with the convention that $\beta^{(p)} = -A\gamma$ if $p = 0$ and that $\beta^{(p)} = 0 \in \mathbb{R}^k$ if $p = P$. In case $p > \max\{p_0(\theta), \mathcal{O}\}$,

$$(48) \qquad G_{\infty,\theta,\sigma,\gamma}(t|p) = \int_{z \le t-\beta^{(p)}} \frac{1 - \Delta_{\sigma\zeta_{\infty,p}}(\nu_p + b_{\infty,p}z, c_p\sigma\xi_{\infty,p})}{1 - \Delta_{\sigma\xi_{\infty,p}}(\nu_p, c_p\sigma\xi_{\infty,p})} \Phi_{\infty,p}(dz),$$

where $\nu_p = \gamma_p + (Q[p\!:\!p]^{-1} Q[p\!:\!\neg p]\gamma[\neg p])_p$. [Note that $\beta^{(p)} = \lim_{n\to\infty} \sqrt{n} \times A(\eta_n(p) - \theta - \gamma/\sqrt{n})$ because $\theta \in M_p$, and that $\nu_p = \lim_{n\to\infty} \sqrt{n}\eta_{n,p}(p)$ in case $\theta \in M_{p-1}$, that is, $p > p_0(\theta)$. Here $\eta_n(p)$ is defined as in (9), but with $\theta + \gamma/\sqrt{n}$ replacing $\theta$.]

If $p > 0$ and if the matrix $A[p]$ has full row rank $k$, then the Lebesgue density $\phi_{\infty,p}(\cdot)$ of $\Phi_{\infty,p}(\cdot)$ exists; the density of (47) is then given by $\phi_{\infty,p}(t - \beta^{(p)})$, while the density of (48) is given by the integrand in (48) times $\phi_{\infty,p}(z)$, evaluated at $z = t - \beta^{(p)}$.

While the limiting c.d.f. in (47) is Gaussian, the limiting c.d.f. in (48) typically is not, an exception being the case where $C_\infty^{(p)} = 0$, that is, when $A\tilde{\theta}(p)$ and $\tilde{\theta}_p(p)$ are asymptotically uncorrelated. In that case, the expressions in (47) and (48) coincide. Also note that the c.d.f. $G_{\infty,\theta,\sigma,\gamma}(t|p)$ has been defined above only for $\theta \in M_p$ (and $\mathcal{O} \le p \le P$). If $\gamma = 0$, we write $G_{\infty,\theta,\sigma}(t|p)$ as shorthand for $G_{\infty,\theta,\sigma,0}(t|p)$ in the following.

Proposition A.1 is restricted to sequences of parameters $\theta + \gamma/\sqrt{n}$ with $p_0(\theta) \le p$. The case where the selected model $M_p$ is an incorrect model for $\theta$, that is, where we have $p_0(\theta) > p$, is analyzed in [10], Proposition 5.1; see also the discussion following Corollary 5.4 in that reference. For the results in the present paper, however, we shall only need to rely on the situation covered by Proposition A.1. The reason essentially is that only over $1/\sqrt{n}$-"neighborhoods" of $M_p$ is the probability of actually selecting the model $M_p$ bounded away from zero. In contrast, for every fixed $\theta \notin M_p$, the probability of selecting the model $M_p$ converges to zero as $n \to \infty$.

PROPOSITION A.2. *Let $p$ satisfy $\mathcal{O} \leq p \leq P$, and let $r_n$ be a sequence of positive real numbers.*

(a) *If $r_n = O(1/\sqrt{n})$ as $n \to \infty$, then*

$$\text{(49)} \qquad \liminf_{n \to \infty} \inf_{\substack{\vartheta \in \mathbf{R}^P \\ \|\vartheta[\neg p]\| < r_n}} P_{n,\vartheta,\sigma}(\hat{p} = p) > 0$$

*holds for every $\sigma$, $0 < \sigma < \infty$. (The infimum in the above display is to be interpreted as extending over $\|\vartheta\| < r_n$ if $p = 0$ and over all of $\mathbf{R}^P$ if $p = P$.) In particular, it follows that $\liminf_{n \to \infty} \inf_{\vartheta \in \mathbf{R}^P, \|\vartheta - \theta\| < r_n} P_{n,\vartheta,\sigma}(\hat{p} = p) > 0$ for each $\theta \in M_p$ and $0 < \sigma < \infty$.*

(b) *Suppose $p < P$ holds. If $\sqrt{n} r_n \to \infty$ as $n \to \infty$, then*

$$\text{(50)} \qquad \lim_{n \to \infty} \inf_{\substack{\vartheta \in \mathbf{R}^P \\ \|\vartheta - \theta\| < r_n}} P_{n,\vartheta,\sigma}(\hat{p} = p) = 0$$

*for each $\theta \in M_p$ and $0 < \sigma < \infty$.*

(c) *If an infimum (resp. supremum) over $\sigma \in [\sigma_*, \sigma^*]$, $0 < \sigma_* \leq \sigma^* < \infty$, is inserted in (49) [resp. (50)] immediately after the $\liminf$ (resp. $\lim$) operator, the result continues to hold.*

PROOF. Let $\vartheta^{(n)}$ be an arbitrary sequence of parameters in $\mathbf{R}^P$. Proposition 5.4 in Leeb [11] together with Remark 5.5 in that reference show that any accumulation point of $P_{n,\vartheta^{(n)},\sigma}(\hat{p} = p)$ is of the form

$$\text{(51)} \qquad (1 - \Delta_{\sigma\xi_{\infty,p}}(v_p, c_p \sigma \xi_{\infty,p})) \prod_{q=p+1}^{P} \Delta_{\sigma\xi_{\infty,q}}(v_q, c_q \sigma \xi_{\infty,q})$$

in case $p > \mathcal{O}$, and of the form

$$\text{(52)} \qquad \prod_{q=p+1}^{P} \Delta_{\sigma\xi_{\infty,q}}(v_q, c_q \sigma \xi_{\infty,q})$$

in case $p = \mathcal{O}$. The quantities $v_q$, $q = p, \ldots, P$, in these displays are accumulation points of $v_q^{(n)} = \sqrt{n}\vartheta_q^{(n)} + \sqrt{n}((X[q]'X[q])^{-1}X[q]'X[\neg q]\vartheta^{(n)}[\neg q])_q$ in $\mathbf{R} \cup \{-\infty, \infty\}$. (In case $q = P$ this expression is to be interpreted as $\sqrt{n}\vartheta_P^{(n)}$ by our conventions.) Observe that the expression in (51) is positive if and only if $|v_q| < \infty$ holds for each $q = p+1, \ldots, P$. The same is true for (52). [In case $p = P$, the expression in (51) is always positive.]

To prove part (a), it suffices to show that any accumulation point of $P_{n,\vartheta^{(n)},\sigma}(\hat{p} = p)$ is positive whenever $\vartheta^{(n)}$ is a sequence satisfying $\|\vartheta^{(n)}[\neg p]\| < r_n$. In case $p = P$ it is easy to see that (51) reduces to $1 - \Delta_{\sigma\xi_{\infty,P}}(v_P, c_P \sigma \xi_{\infty,P})$, which is bounded from below by the positive constant $1 - \Delta_{\sigma\xi_{\infty,P}}(0, c_P \sigma \xi_{\infty,P})$. In case $p < P$ note that $\sqrt{n}\vartheta^{(n)}[\neg p]$ is a bounded sequence and, hence, $v_q^{(n)}$ is bounded for each $q = p+1, \ldots, P$. It follows that $|v_q| < \infty$ holds for $q = p+1, \ldots, P$. This completes the proof of part (a).

To prove part (b), let $\vartheta^{(n)}$ be given by $\vartheta^{(n)}[P-1] = \theta[P-1]$ and $\vartheta_P^{(n)} = r_n/2$. Clearly, then $\|\vartheta^{(n)} - \theta\| < r_n$ is satisfied. Moreover, $\sqrt{n}\vartheta_P^{(n)} = \sqrt{n}r_n/2$ converges to $v_P = \infty$. It follows that $\lim_{n \to \infty} P_{n,\vartheta^{(n)},\sigma}(\hat{p} = p) = 0$ whenever $p < P$.

Part (c) is proved analogously. $\square$

## APPENDIX B: PROOFS FOR SECTION 2.2.1

In the proofs below it will be convenient to show the dependence of $\Phi_{n,p}(t)$ and $\Phi_{\infty,p}(t)$ on $\sigma$ in the notation. Thus, in the following we shall write $\Phi_{n,p,\sigma}(t)$ and $\Phi_{\infty,p,\sigma}(t)$ for the c.d.f. of a mean zero $k$-variate Gaussian random vector with variance–covariance matrix $\sigma^2 A[p](X[p]'X[p]/n)^{-1}A[p]'$ and $\sigma^2 A[p]Q[p:p]^{-1}A[p]'$, respectively. For convenience, let $\Phi_{n,0,\sigma}(t)$ and $\Phi_{\infty,0,\sigma}(t)$ denote c.d.f.s of a point-mass at zero in $\mathbf{R}^k$. The following lemma is elementary to prove, if we observe that in case $\text{rank}(A[p]) = k$ the convergence $b_{n,p} \to b_{\infty,p}$ holds, since the generalized inverses in the definitions of these quantities then reduce to the usual inverse.

LEMMA B.1. *Suppose $p > \mathcal{O}$ and that $\text{rank}(A[p]) = k$. Define $S_{n,p}(z, \sigma) = \frac{1 - \Delta_{\sigma\zeta_{n,p}}(b_{n,p}z, c_p\sigma\xi_{n,p})}{1 - \Delta_{\sigma\xi_{n,p}}(0, c_p\sigma\xi_{n,p})}$ and $S_{\infty,p}(z, \sigma) = \frac{1 - \Delta_{\sigma\zeta_{\infty,p}}(b_{\infty,p}z, c_p\sigma\xi_{\infty,p})}{1 - \Delta_{\sigma\xi_{\infty,p}}(0, c_p\sigma\xi_{\infty,p})}$ for $z \in \mathbf{R}^k$, $0 < \sigma < \infty$. Let $\sigma^{(n)}$ converge to $\sigma$, $0 < \sigma < \infty$. Then $S_{n,p}(z, \sigma^{(n)})$ converges to $S_{\infty,p}(z, \sigma)$ for every $z \in \mathbf{R}^k$ if $\zeta_{\infty,p} \neq 0$, and for every $z \in \mathbf{R}^k$ except possibly for $z$ satisfying $|b_{\infty,p}z| = c_p\sigma\xi_{\infty,p}$ if $\zeta_{\infty,p} = 0$. (The exceptional set has Lebesgue measure zero since $c_p\sigma\xi_{\infty,p} > 0$.)*

LEMMA B.2. *Let $(\Omega, \mathcal{A})$ and $(\Xi, \mathcal{B})$ be measurable spaces and let $\Psi: \Omega \to \Xi$ be a measurable function. Suppose $\mu_n$ and $\mu$ are probability measures on $(\Omega, \mathcal{A})$ satisfying $\|\mu_n - \mu\|_{\text{TV}} \to 0$. Let $\rho_n$ be the probability measure induced by $\mu_n$ and $\Psi$, that is, $\rho_n(B) = \mu_n(\Psi^{-1}(B))$ for $B \in \mathcal{B}$. Then $\rho_n$ converges to a probability measure $\rho$ with respect to the total variation distance and $\rho$ is the measure induced by $\mu$ and $\Psi$.*

Lemma B.2 follows immediately from $\|\rho_n - \rho\|_{\text{TV}} \leq \|\mu_n - \mu\|_{\text{TV}}$. The following observation is useful in the proof of Proposition 2.1 below: Since the proposition depends on $Y$ only through its distribution (cf. Remark 4.6), we may assume without loss of generality that the errors in (5) are given by $u_t = \sigma\varepsilon_t$, $t \in \mathbf{N}$, with i.i.d. $\varepsilon_t$ that are standard normal. In particular, all random variables involved are then defined on the same probability space.

PROOF OF PROPOSITION 2.1. We consider first the case $p > \mathcal{O}$ and assume for the moment that the matrix $A[p]$ has full row rank $k$. Then $\Phi_{n,p,\sigma}(\cdot)$ and $\Phi_{\infty,p,\sigma}(\cdot)$ possess densities $\phi_{n,p,\sigma}(\cdot)$ and $\phi_{\infty,p,\sigma}(\cdot)$, respectively, with respect to Lebesgue measure on $\mathbf{R}^k$. Since $\hat{\sigma} \to \sigma$ in $P_{n,\theta,\sigma}$-probability, each subsequence contains a further subsequence along which $\hat{\sigma} \to \sigma$ almost surely (with respect to the probability measure on the common probability space supporting all random variables involved), and we restrict ourselves to this further subsequence for the moment. In particular, we write $\{\hat{\sigma} \to \sigma\}$ for the event that $\hat{\sigma}$ converges to $\sigma$ along the subsequence under consideration; clearly, the event $\{\hat{\sigma} \to \sigma\}$ has probability one. Also note that we can assume without loss of generality that $\hat{\sigma} > 0$ holds on this event (at least from some data-dependent $n$ onward), since $\sigma > 0$ holds. Lemma B.1 now shows that on the event $\{\hat{\sigma} \to \sigma\}$ the function $S_{n,p}(z, \hat{\sigma})\phi_{n,p,\hat{\sigma}}(z)$ converges to $S_{\infty,p}(z, \sigma)\phi_{\infty,p,\sigma}(z)$ for every $z$ except for a set of Lebesgue measure zero. Observe that both functions are probability densities with respect to Lebesgue measure on $\mathbf{R}^k$; see the discussion prior to Proposition 2.1. In view of Scheffé's lemma, they hence converge in absolute mean. By the same argument, $\phi_{n,p,\hat{\sigma}}(\cdot)$ also converges to $\phi_{\infty,p,\sigma}(\cdot)$ in absolute mean. Note that the absolute mean convergence of the densities translates into convergence in total variation for the corresponding c.d.f.s. Now $\check{G}_n(t|p) = \Phi_{n,p,\hat{\sigma}}(t)$ in case the auxiliary procedure decides for $p_0(\theta) = p$, and $\check{G}_n(t|p) = \int_{z \in \mathbf{R}^k, z \leq t} S_{n,p}(z, \hat{\sigma})\phi_{n,p,\hat{\sigma}}(z) \, dz$ otherwise. Since the auxiliary procedure decides consistently between $p_0(\theta) = p$ and $p_0(\theta) < p$ for every $\theta \in M_p$, it follows that (18) holds along the subsequence under consideration in case $p > \mathcal{O}$ and if $A[p]$ has rank $k$. Of course, this already proves (18) in case $p > \mathcal{O}$ and $A[p]$ has rank $k$.

In case $p > \mathcal{O}$ but where the matrix $A[p]$ does not have full row rank $k$, let $G^I_{n,\theta,\sigma}(t|p)$, $\check{G}^I_n(t|p)$ and $G^I_{\infty,\theta,\sigma}(t|p)$ be defined in exactly the same way as $G_{n,\theta,\sigma}(t|p)$, $\check{G}_n(t|p)$ and $G_{\infty,\theta,\sigma}(t|p)$, respectively, except that the $p \times P$ matrix $(I_p : 0)$ replaces $A$. Note that then $I_p$ replaces $A[p]$ (and that the value of $k$ changes to $p$). Since the matrix $I_p$ has full row rank $p$, the preceding paragraph shows that (18) holds with $\check{G}^I_n(t|p)$ and $G^I_{\infty,\theta,\sigma}(t|p)$ replacing $\check{G}_n(t|p)$ and $G_{\infty,\theta,\sigma}(t|p)$, respectively. But $\check{G}_n(t|p)$ and $G_{\infty,\theta,\sigma}(t|p)$, respectively, are the c.d.f.s of the image measures of $\check{G}^I_n(t|p)$ and $G^I_{\infty,\theta,\sigma}(t|p)$ induced by the linear mapping $x \mapsto A[p]x$, $x \in \mathbf{R}^p$. [This is obvious for $\check{G}_n(t|p)$ because of its interpretation as the conditional c.d.f. $G^*_{n,\theta,\hat{\sigma}}(t|p)$ in (13) of [10] if $\hat{\sigma} > 0$; it is trivial if $\hat{\sigma} = 0$. Observe further that $G_{n,\theta,\sigma}(t|p)$ is clearly the c.d.f. of the induced measure obtained from the c.d.f. $G^I_{n,\theta,\sigma}(t|p)$. Since $G^I_{n,\theta,\sigma}(t|p) \to G^I_{\infty,\theta,\sigma}(t|p)$ and $G_{n,\theta,\sigma}(t|p) \to G_{\infty,\theta,\sigma}(t|p)$ with respect to total variation distance for $\theta \in M_p$ by Proposition A.1, an application of Lemma B.2 shows that $G_{\infty,\theta,\sigma}(t|p)$ is indeed the c.d.f. of the induced measure obtained from the c.d.f. $G^I_{\infty,\theta,\sigma}(t|p)$.] Therefore, the total variation distance of $\check{G}_n(t|p)$ and $G_{\infty,\theta,\sigma}(t|p)$ is bounded from above by that of $\check{G}^I_n(t|p)$ and $G^I_{\infty,\theta,\sigma}(t|p)$. This proves (18) also in this case.

In the case $p = \mathcal{O} > 0$, note that $G_{\infty,\theta,\sigma}(t|p)$ is given by (47) for $\theta \in M_p$. The result in (18) then follows in a similar way, observing that in case $A[p]$ has full row rank $k$ (again after passing to appropriate subsequences), $\phi_{n,p,\hat{\sigma}}(\cdot)$ converges to $\phi_{\infty,p,\sigma}(\cdot)$ in absolute mean on the event $\{\hat{\sigma} \to \sigma\}$ as defined above. The case where $p = \mathcal{O} = 0$ is trivial, because both c.d.f.s in (18) coincide and are equal to the c.d.f. of point-mass at zero in $\mathbf{R}^k$. This completes the proof of (18).

The validity of (17) now follows for $\theta \in M_p$ since $G_{\infty,\theta,\sigma}(t|p)$ is then the limit of $G_{n,\theta,\sigma}(t|p)$ with respect to the total variation distance; see Proposition A.1. Finally, the claim regarding "conditional consistency" follows from (17) and (18) in view of Proposition A.2(a). □

PROOF OF COROLLARY 2.2. Observe that

$$P_{n,\theta,\sigma}(\|\check{G}_n(\cdot|\hat{p}) - G_{n,\theta,\sigma}(\cdot|\hat{p})\|_{\mathrm{TV}} > \delta)$$

$$= \sum_{p=\mathcal{O}}^{P} P_{n,\theta,\sigma}(\|\check{G}_n(\cdot|p) - G_{n,\theta,\sigma}(\cdot|p)\|_{\mathrm{TV}} > \delta, \hat{p} = p)$$

$$\leq \sum_{\mathcal{O} \leq p < p_0(\theta)} P_{n,\theta,\sigma}(\hat{p} = p)$$

$$+ \sum_{p \geq p_0(\theta)} P_{n,\theta,\sigma}(\|\check{G}_n(\cdot|p) - G_{n,\theta,\sigma}(\cdot|p)\|_{\mathrm{TV}} > \delta).$$

Each term in the first sum on the far right-hand side of the above display now obviously converges to zero (cf. [11], Corollary 5.6 and (5.7)), whereas every term in the second sum converges to zero by Proposition 2.1. □

## APPENDIX C: PROOFS FOR SECTIONS 2.2.2 AND 4.1

Since the results in Section 2.2.2 rely on those in Section 4.1, the latter ones are proved first. Some of the proofs rely on auxiliary results collected in Appendix D. We start with some preparatory remarks. The total variation distance between $P_{n,\theta,\sigma}$ and $P_{n,\vartheta,\sigma}$ satisfies $\|P_{n,\theta,\sigma} - P_{n,\vartheta,\sigma}\|_{\mathrm{TV}} \leq 2\Phi(\|\theta - \vartheta\|\lambda_{\max}^{1/2}(X'X)/2\sigma) - 1$; furthermore, if $\theta^{(n)}$ and $\vartheta^{(n)}$ satisfy $\|\theta^{(n)} - \vartheta^{(n)}\| = O(n^{-1/2})$, the sequence $P_{n,\vartheta^{(n)},\sigma}$ is contiguous with respect to the sequence $P_{n,\theta^{(n)},\sigma}$. This follows exactly in the same way as Lemma A.1 in [16]. We also need the following lemma.

LEMMA C.1. *Let $p$ satisfy $\mathcal{O} < p \leq P$. Suppose $\theta \in M_{p-1}$, $0 < \sigma < \infty$ and $0 < \rho < \infty$. Then*

$$\liminf_{n \to \infty} \inf_{\substack{\vartheta \in M_p \\ \|\vartheta - \theta\| < \rho/\sqrt{n}}} P_{n,\vartheta,\sigma}(\hat{p} = p)$$

$$= (1 - \Delta_{\sigma\xi_{\infty,p}}(0, c_p\sigma\xi_{\infty,p})) \prod_{q=p+1}^{P} \Delta_{\sigma\xi_{\infty,q}}(0, c_q\sigma\xi_{\infty,q})$$

(53)

$$= 2(1 - \Phi(c_p)) \prod_{q=p+1}^{P} (2\Phi(c_q) - 1)$$

$$= \lim_{n \to \infty} P_{n,\theta,\sigma}(\hat{p} = p) > 0.$$

PROOF. We proceed similarly as in the proof of Proposition A.2, observing that now the quantities $v_q$, $q > p$, are all equal to zero since $\vartheta^{(n)} \in M_p$. Since $(1 - \Delta_{\sigma\xi_{\infty,p}}(v_p, c_p\sigma\xi_{\infty,p}))$ is minimal for $v_p = 0$, we see that the right-hand side of (53), which obviously is positive, is a lower bound for the left-hand side. Using (5.7) in [11] and observing that $\theta \in M_{p-1}$ completes the proof. $\square$

PROOF OF THEOREM 4.1. We first prove (33) and (34). For this purpose, we make use of Lemma 3.1 in [16] with $\alpha = \theta \in M_{p-1}$, $B = M_p$, $B_n = \{\vartheta \in M_p : \|\vartheta - \theta\| < \rho_0 n^{-1/2}\}$, $\beta = \vartheta$, $\varphi_n(\beta) = G_{n,\vartheta,\sigma}(t|p)$ and $\hat{\varphi}_n = \hat{G}_n(t|p)$, where $\rho_0$, $0 < \rho_0 < \infty$, will be chosen shortly (and $\sigma$ is held fixed). The contiguity assumption of this lemma is satisfied in view of the preparatory remark above. It hence only remains to show that there exists a value of $\rho_0$, $0 < \rho_0 < \infty$, such that $\delta^*$ in Lemma 3.1 of [16] [which represents the limit inferior of the oscillation of $\varphi_n(\cdot)$ over $B_n$] is positive. Applying Lemma 3.5(a) of [16] with $\zeta_n = \rho_0 n^{-1/2}$ and the set $G_0$ equal to the set $G$, it remains, in light of Proposition A.1, to show that there exists a $\rho_0$, $0 < \rho_0 < \infty$, such that $G_{\infty,\theta,\sigma,\gamma}(t|p)$ as a function of $\gamma$ is nonconstant on the set $\{\gamma \in M_p : \|\gamma\| < \rho_0\}$. In view of Lemma 3.1 of [16], the corresponding $\delta_0$ can then be chosen as any positive number less than one-half of the oscillation of $G_{\infty,\theta,\sigma,\gamma}(t|p)$ over this set. That such a $\rho_0$ indeed exists follows from Lemma D.1 in Appendix D. Furthermore, observe that $G_{\infty,\theta,\sigma,\cdot}(t|p)$ is given by (48) for $\theta \in M_{p-1}$ and, hence, does not depend on $\theta$, but only on $t, Q, A, \sigma$ and $c_p$. As a consequence, $\rho_0$ and $\delta_0$ can be chosen such that they also depend only on these quantities. This completes the proof of (33) and (34).

To prove (35), we use Corollary 3.4 in [16] with the same identification of notation as above, with $\zeta_n = \rho_0 n^{-1/2}$, and with $V = M_p$ (viewed as a vector space isomorphic to $\mathbf{R}^p$). The asymptotic uniform equicontinuity condition in that corollary is then satisfied in view of $\|P_{n,\theta,\sigma} - P_{n,\vartheta,\sigma}\|_{\mathrm{TV}} \leq 2\Phi(\|\theta - \vartheta\|\lambda_{\max}^{1/2}(X'X)/2\sigma) - 1$. Applying Corollary 3.4 in [16] then establishes (35). This completes the proof of part (a).

Part (b) is proved exactly as part (a), making additional use of Corollary C.2 and Remark C.1 in [16]. The events $E_n$ appearing in this corollary are given here by $\{\hat{p}=p\}$. Clearly, $P_{n,\vartheta,\sigma}(\hat{p}=p)$ is always positive. The constant $M$ in Corollary C.2 of [16] is now given by the right-hand side of (53) above. $\square$

PROOF OF THEOREM 4.2. We again use results from [16], this time with the identification $\alpha = \theta \in M_p$, $B = M_{q^*}$, $B_n = \{\vartheta \in M_{q^*} : \|\vartheta - \theta\| < \rho_0 n^{-1/2}\}$, $\beta = \vartheta$, $\varphi_n(\beta) = G_{n,\vartheta,\sigma}(t|p)$, $\widehat{\varphi}_n = \hat{G}_n(t|p)$, $V = M_{q^*}$ and $\zeta_n = \rho_0 n^{-1/2}$ (again $\sigma$ is held fixed). The proof of part (a) is then similar to the proof of part (a) of Theorem 4.1, except for using Lemma D.2 instead of Lemma D.1 and except for the fact that the argument that $\rho_0$ and $\delta_0$ only depend on $t, Q, A, \sigma$ and $c_p$ is now slightly more complex, since $G_{\infty,\theta,\sigma,\cdot}(t|p)$ for $\theta \in M_p$ depends on $\theta$. However, observe that $G_{\infty,\theta,\sigma,\cdot}(t|p)$ as a function of $\theta \in M_p$ can follow only two different formulae which themselves do not depend on $\theta$; see (47) and (48).

Part (b) is proved exactly as the corresponding part of Theorem 4.1, except that positivity of the constant $M = \liminf_{n\to\infty} \inf_{\vartheta \in M_{q^*}, \|\vartheta-\theta\|<\rho/\sqrt{n}} P_{n,\vartheta,\sigma}(\hat{p}=p)$ follows now since $M$ is bounded from below by the expression in part (a) of Proposition A.2. $\square$

PROOF OF PROPOSITION 4.3. See [14]. $\square$

PROOF OF PROPOSITION 4.4. That part (a) implies part (b) follows from (20) in [10], observing that $C_n^{(q)} \to C_\infty^{(q)}$ and that $\eta_{n,q}(q)$ converges to a finite limit. The reverse implication follows by passing to the limit in (20) of [10] and observing that, by suitable choice of $\theta \in \mathbf{R}^P$, the limit of $(\eta_{n,p+1}(p+1), \ldots, \eta_{n,P}(P))'$ can take on the value of every standard basis vector in $\mathbf{R}^{P-p}$. To prove the equivalence of parts (a) and (c), we use Proposition 3.1 in [10] and equation (19) in that paper to obtain $\sum_{r=1}^{q} \sigma^2 \xi_{\infty,r}^{-2} C_\infty^{(r)} C_\infty^{(r)\prime}$ as the formula for the asymptotic variance–covariance matrix of $\sqrt{n}A\tilde{\theta}(q)$. Since the terms in this sum are nonnegative definite, the equivalence follows. The final claims regarding the cases $p = P$ and $p = 0$ are either obvious or follow immediately from the representation of the asymptotic variance–covariance matrix of $\sqrt{n}A\tilde{\theta}(q)$ just given. $\square$

PROOF OF THEOREM 2.3. In view of the definition of $G_{n,\vartheta,\sigma}(t|\hat{p})$, we have

$$|\hat{G}_n(t|\hat{p}) - G_{n,\vartheta,\sigma}(t|\hat{p})| = \sum_{p=\mathcal{O}}^{P} |\hat{G}_n(t|\hat{p}) - G_{n,\vartheta,\sigma}(t|p)|\mathbf{1}(\hat{p}=p)$$

$$\geq |\hat{G}_n(t|\hat{p}) - G_{n,\vartheta,\sigma}(t|q^*)|\mathbf{1}(\hat{p}=q^*).$$

Hence, for every $\vartheta \in \mathbf{R}^P$ and every $\delta > 0$,

$$P_{n,\vartheta,\sigma}(|\hat{G}_n(t|\hat{p}) - G_{n,\vartheta,\sigma}(t|\hat{p})| > \delta)$$
$$\geq P_{n,\vartheta,\sigma}(|\hat{G}_n(t|\hat{p}) - G_{n,\vartheta,\sigma}(t|q^*)| > \delta|\hat{p} = q^*) P_{n,\vartheta,\sigma}(\hat{p} = q^*),$$

observing that the conditional probabilities are well defined since $P_{n,\vartheta,\sigma}(\hat{p} = q^*)$ is always positive (cf. [11], Section 3.2). This implies

$$\liminf_{n \to \infty} \sup_{\substack{\vartheta \in M_{q^*} \\ \|\vartheta - \theta\| < \rho_0/\sqrt{n}}} P_{n,\vartheta,\sigma}(|\hat{G}_n(t|\hat{p}) - G_{n,\vartheta,\sigma}(t|\hat{p})| > \delta)$$

$$(54) \quad \geq \left[ \liminf_{n \to \infty} \sup_{\substack{\vartheta \in M_{q^*} \\ \|\vartheta - \theta\| < \rho_0/\sqrt{n}}} P_{n,\vartheta,\sigma}(|\hat{G}_n(t|\hat{p}) - G_{n,\vartheta,\sigma}(t|q^*)| > \delta|\hat{p} = q^*) \right.$$

$$\left. \times \liminf_{n \to \infty} \inf_{\substack{\vartheta \in M_{q^*} \\ \|\vartheta - \theta\| < \rho_0/\sqrt{n}}} P_{n,\vartheta,\sigma}(\hat{p} = q^*) \right].$$

Lemma C.1 above shows that

$$\liminf_{n \to \infty} \inf_{\substack{\vartheta \in M_{q^*} \\ \|\vartheta - \theta\| < \rho_0/\sqrt{n}}} P_{n,\vartheta,\sigma}(\hat{p} = q^*)$$

$$= 2(1 - \Phi(c_{q^*})) \prod_{q=q^*+1}^{P} (2\Phi(c_q) - 1) = \lim_{n \to \infty} P_{n,\theta,\sigma}(\hat{p} = q^*),$$

which obviously is positive. Suppose now that $\hat{G}_n(t|\hat{p})$ satisfies (19). Then it also satisfies $P_{n,\theta,\sigma}(|\hat{G}_n(t|\hat{p}) - G_{n,\theta,\sigma}(t|q^*)| > \delta|\hat{p} = q^*) \xrightarrow{n \to \infty} 0$, since the probability $P_{n,\theta,\sigma}(\hat{p} = q^*)$ of the conditioning event is bounded away from zero as just shown. Since $q^* > \mathcal{O}$ is the maximal model order $q$ with the property that $C_\infty^{(q)} \neq 0$, the condition on $t$ in Theorem 4.1 is satisfied for every $t \in \mathbf{R}^k$. Hence, we may apply Theorem 4.1(b) with $p = q^*$ to the first term in the product on the right-hand side of (54) since $\hat{G}_n(t|\hat{p})$ can certainly also be viewed as an estimator of $G_{n,\theta,\sigma}(t|q^*)$. This establishes (20) with the same $\delta_0$ and $\rho_0$ as in Theorem 4.1(b). Furthermore, note that (54) remains valid if an infimum extending over all estimators is inserted between the limit inferior and the supremum on both sides of (54). Again applying Theorem 4.1(b) with $p = q^*$ completes the proof of (21)–(22). □

PROOF OF PROPOSITION 2.4. See [14]. □

## APPENDIX D: AUXILIARY LEMMATA FOR APPENDIX C

LEMMA D.1. *Let $p$ satisfy $\mathcal{O} < p \leq P$, and assume that $A\tilde{\theta}(p)$ and $\tilde{\theta}_p(p)$ are asymptotically correlated, that is, $C_\infty^{(p)} \neq 0$. Moreover, let $\theta \in M_{p-1}$, let $\sigma$ satisfy $0 < \sigma < \infty$ and let $t \in \mathbf{R}^k$ be such that the set $\{z \in \mathbf{R}^p : A[p]z \leq t\}$ has positive Lebesgue measure in $\mathbf{R}^p$ (which is satisfied for all $t \in \mathbf{R}^k$ if, e.g., $\mathrm{rank}(A[p]) = k$). Then $G_{\infty,\theta,\sigma,\gamma}(t|p)$ is nonconstant as a function of $\gamma \in M_p$.*

LEMMA D.2. *Let $p$ satisfy $\mathcal{O} \leq p < P$, assume that $A\tilde{\theta}(q)$ and $\tilde{\theta}_q(q)$ are asymptotically correlated, that is, $C_\infty^{(q)} \neq 0$, for some $q$ satisfying $p < q \leq P$, and let $q^*$ denote the largest $q$ with this property. Moreover, let $t \in \mathbf{R}^k$, let $\theta \in M_p$ and let $\sigma$ satisfy $0 < \sigma < \infty$. Then $G_{\infty,\theta,\sigma,\gamma}(t|p)$ is nonconstant as a function of $\gamma \in M_{q^*}$.*

Before we prove the above lemmata, we provide a representation of $G_{\infty,\theta,\sigma,\gamma}(t|p)$ for $p > 0$ that will be useful in the following. For $0 < p \leq P$, define $Z_p = \sum_{r=1}^p \xi_{\infty,r}^{-2} C_\infty^{(r)} W_r$, where $C_\infty^{(r)}$ has been defined after (45) and the random variables $W_r$ are independent and normally distributed with mean zero and variances $\sigma^2 \xi_{\infty,r}^2$. For convenience, let $Z_0$ denote the zero vector in $\mathbf{R}^k$. Observe that $Z_p$, $p > 0$, is normally distributed with mean zero and variance–covariance matrix $\sigma^2 A[p]Q[p:p]^{-1}A[p]'$, since we have shown in the proof of Proposition 4.4 that the asymptotic variance–covariance matrix of $\sqrt{n}A\tilde{\theta}(p)$ can be expressed as $\sum_{r=1}^p \sigma^2 \xi_{\infty,r}^{-2} C_\infty^{(r)} C_\infty^{(r)'}$. Also, the joint distribution of $Z_p$ and the set of variables $W_r$, $1 \leq r \leq P$, is normal, with the covariance vector between $Z_p$ and $W_r$ given by $\sigma^2 C_\infty^{(r)}$ in case $r \leq p$; otherwise $Z_p$ and $W_r$ are independent. Define the constants $\nu_r = \gamma_r + (Q[r:r]^{-1}Q[r:\neg r]\gamma[\neg r])_r$ for $0 < r \leq P$. It is now easy to see that $\beta^{(p)}$ defined in Proposition A.1 equals $-\sum_{r=p+1}^P \xi_{\infty,r}^{-2} C_\infty^{(r)} \nu_r$. [This is seen as follows: It was noted in Proposition A.1 that $\beta^{(p)} = \lim_{n \to \infty} \sqrt{n}A(\eta_n(p) - \theta - \gamma/\sqrt{n})$ for $\theta \in M_p$, when $\eta_n(p)$ is defined as in (9), but with $\theta + \gamma/\sqrt{n}$ replacing $\theta$. Using the representation (20) of [10] and taking limits, the result follows if we observe that $\sqrt{n}\eta_{n,r}(r) \longrightarrow \nu_r$ for $\theta \in M_p$ and $r > p$.] In view of (47), the c.d.f. $G_{\infty,\theta,\sigma,\gamma}(t|p)$ can now equivalently be written as

$$(55) \quad G_{\infty,\theta,\sigma,\gamma}(t|p) = P\left(Z_p \leq t + \sum_{r=p+1}^P \xi_{\infty,r}^{-2} C_\infty^{(r)} \nu_r\right)$$

in case $p = \max\{p_0(\theta), \mathcal{O}\} > 0$, and (55) trivially holds in case $p = 0$. In case $p > \max\{p_0(\theta), \mathcal{O}\}$ the c.d.f. $G_{\infty,\theta,\sigma,\gamma}(t|p)$ is given by (48), and it is elementary but tedious to show, following the steps in Section 3.1 of [10], that this is equivalent to

$$(56) \quad G_{\infty,\theta,\sigma,\gamma}(t|p) = P\left(Z_p \leq t + \sum_{r=p+1}^P \xi_{\infty,r}^{-2} C_\infty^{(r)} \nu_r \,\bigg|\, |W_p + \nu_p| \geq c_p \sigma \xi_{\infty,p}\right).$$

(This can also be derived from the fact that the distribution of $(Z_p', W_p + \nu_p, \ldots, W_P + \nu_P)'$ represents the limiting distribution of $\sqrt{n}(A(\tilde{\theta}(p) - \eta_n(p))', \tilde{\theta}_p(p), \ldots, \tilde{\theta}_P(P))'$ under $P_{n, \theta + \gamma/\sqrt{n}, \sigma}$ with $\theta \in M_{p-1}$ [and $\eta_n(p)$ defined as in (9), but with $\theta + \gamma/\sqrt{n}$ replacing $\theta$].)

PROOF OF LEMMA D.1. Since $\theta \in M_{p-1}$, $G_{\infty,\theta,\sigma,\gamma}(t|p)$ is given by (56). For $\gamma \in M_p$, the quantities $\nu_{p+1}, \ldots, \nu_P$ are easily seen to be zero, while $\nu_p$ equals $\gamma_p$. This leads to

$$G_{\infty,\theta,\sigma,\gamma}(t|p) = P(Z_p \leq t ||W_p + \gamma_p| \geq c_p \sigma \xi_{\infty,p})$$

for $\gamma \in M_p$. Since $Z_p = Z_{p-1} + \xi_{\infty,p}^{-2} C_\infty^{(p)} W_p$, we obtain

(57) $\quad G_{\infty,\theta,\sigma,\gamma}(t|p) = P(Z_{p-1} + \xi_{\infty,p}^{-2} C_\infty^{(p)} W_p \leq t ||W_p + \gamma_p| \geq c_p \sigma \xi_{\infty,p}).$

Assume now that (57) is constant in $\gamma_p \in \mathbf{R}$. Using Lemma D.3 below with $Z_{p-1} - t$, $W_p$, $-\xi_{\infty,p}^{-2} C_\infty^{(p)}$, $-\gamma_p$ and $c_p \sigma \xi_{\infty,p}$ replacing $Z$, $W$, $C$, $x$ and $\delta$, respectively, we obtain that either $P(Z_p \leq t) = 0$ or that $\xi_{\infty,p}^{-2} C_\infty^{(p)} = 0$. By assumption of the lemma, the set $\{z \in \mathbf{R}^p : A[p]z \leq t\}$ has positive Lebesgue measure. Hence, $P(Z_p \leq t)$ must be positive. (To see why, note that $Z_p$ is concentrated in the column space of $A[p]$, and that $Z_p$ is nondegenerate *within* the column-space of $A[p]$.) It would follow that $\xi_{\infty,p}^{-2} C_\infty^{(p)} = 0$, contradicting the assumption that $A\tilde{\theta}(p)$ and $\tilde{\theta}_p(p)$ are asymptotically correlated. □

PROOF OF LEMMA D.2. By the assumptions on $q^*$, note that either $q^* = P$ or that $C_\infty^{(r)} = 0$ for each $r = q^* + 1, \ldots, P$. Consider first the case $p = \max\{p_0(\theta), \mathcal{O}\}$. By (55), we have $G_{\infty,\theta,\sigma,\gamma}(t|p) = P(Z_p \leq t + \sum_{r=p+1}^{q^*} \xi_{\infty,r}^{-2} \times C_\infty^{(r)} \nu_r)$. Observe that $(\nu_{p+1}, \ldots, \nu_{q^*})'$ varies in all of $\mathbf{R}^{q^* - p}$ when $\gamma$ varies in $M_{q^*}$. Hence, the last mentioned probability goes to zero along an appropriate sequence of $(\nu_{p+1}, \ldots, \nu_{q^*})'$ (viz., a sequence along which at least one coordinate of $t + \sum_{r=p+1}^{q^*} \xi_{\infty,r}^{-2} C_\infty^{(r)} \nu_r$ goes to $-\infty$). Since $Z_{q^*} = Z_p + \sum_{r=p+1}^{q^*} \xi_{\infty,r}^{-2} C_\infty^{(r)} W_r$ and since the $W_r$, $r = p + 1, \ldots, P$, are independent of $Z_p$, the c.d.f. $G_{\infty,\theta,\sigma,\gamma}(t|p)$ can also be written as a (regular) conditional probability

(58) $\quad G_{\infty,\theta,\sigma,\gamma}(t|p) = P(Z_{q^*} \leq t | W_{p+1} = -\nu_{p+1}, \ldots, W_{q^*} = -\nu_{q^*}).$

Suppose now that $G_{\infty,\theta,\sigma,\gamma}(t|p)$ is constant in $\gamma \in M_{q^*}$, or equivalently, is constant when $(\nu_{p+1}, \ldots, \nu_{q^*})'$ varies in all of $\mathbf{R}^{q^* - p}$. It follows from the above discussion that the conditional probability in (58) is then zero for all $(\nu_{p+1}, \ldots, \nu_{q^*})' \in \mathbf{R}^{q^* - p}$. By integration with respect to the distribution of $(W_{p+1}, \ldots, W_{q^*})$, we obtain that $P(Z_{q^*} \leq t) = 0$. From Proposition 4.4(c), it follows that $Z_{q^*}$ has a nonsingular normal distribution on $\mathbf{R}^k$, which contradicts $P(Z_{q^*} \leq t) = 0$. This proves the lemma in case $p = \max\{p_0(\theta), \mathcal{O}\}$. Consider next the case $p > \max\{p_0(\theta), \mathcal{O}\}$ and assume that $G_{\infty, \theta, \sigma, \gamma}(t|p)$ is constant in $\gamma \in M_{q^*}$. Now $G_{\infty, \theta, \sigma, \gamma}(t|p)$ is given by (56). Letting $\gamma_p \to \infty$, $\nu_p$ converges to $\infty$ as well, and the expression in (56) converges to that in (55). Hence, (55) would have to be constant as a function of $(\nu_{p+1}, \ldots, \nu_{q^*})'$ (note that $(\nu_{p+1}, \ldots, \nu_{q^*})'$ depends only on $\gamma[\neg p]$ but not on $\gamma_p$), which already has been shown to lead to a contradiction. $\square$

LEMMA D.3. *Let $Z$ be a random vector with values in $\mathbf{R}^k$, let $W$ be a univariate random variable independent of $Z$ and assume that $W$ has a Lebesgue density which is positive almost everywhere. Furthermore, let $C \in \mathbf{R}^k$ and let $\delta > 0$. Then $P(Z \leq CW||W - x| \geq \delta)$ is constant in $x \in \mathbf{R}$ if and only if $P(Z \leq CW) = 0$ or $C = 0$.*

PROOF. If $C = 0$, then $P(Z \leq CW||W - x| \geq \delta)$ equals $P(Z \leq 0)$, which is constant in $x$. If $P(Z \leq CW) = 0$, obviously also $P(Z \leq CW||W - x| \geq \delta) = 0$, and hence, is constant in $x$. Conversely, assume that $P(Z \leq CW||W - x| \geq \delta) = P(Z \leq CW||W - x'| \geq \delta)$ for each $x, x' \in \mathbf{R}$. Letting $x' \to \infty$ implies that

$$\frac{P(Z \leq CW, |W - x| \geq \delta)}{P(|W - x| \geq \delta)} = P(Z \leq CW)$$

holds for each $x \in \mathbf{R}$. This is equivalent to

(59) $\qquad P(Z \leq CW, W \in B) = P(Z \leq CW)P(W \in B),$

for all sets $B$ of the form $B = (x - \delta, x + \delta)$ with $x \in \mathbf{R}$. Since both sides in (59) are sigma-additive set functions and since $W$ is absolutely continuous with respect to Lebesgue measure, both set functions also agree on all sets of the form $(-\infty, x + \delta]$, and hence, on the entire Borel sigma-field on $\mathbf{R}$. This implies independence of $\{Z \leq CW\}$ and $W$. In particular, we have

$$P(Z \leq CW) = P(Z \leq CW|W = w)$$

for almost all $w \in \mathbf{R}$. Furthermore, by the assumed independence of $Z$ and $W$, we have

$$P(Z \leq CW) = P(Z \leq CW|W = w) = P(Z \leq Cw)$$

for almost all $w \in \mathbf{R}$. Now if $C \neq 0$, the right-hand side of the above display goes to zero either for $w \to \infty$ or for $w \to -\infty$, implying that $P(Z \leq CW) = 0$. $\square$

## APPENDIX E: PROOFS FOR SECTIONS 3 AND 4.2

PROOF OF THEOREM 4.5. After rearranging the elements of $\theta$ (and hence, the regressors) and correspondingly rearranging the rows of the matrix $A$ if necessary, we may assume without loss of generality that $\mathfrak{r}_* = (1,\ldots,1,0)$, and hence, that $i(\mathfrak{r}_*) = P$. That is, $M_{\mathfrak{r}_*} = M_{P-1}$ and $M_{\mathfrak{r}_{\text{full}}} = M_P$. Furthermore, note that after this arrangement $C_\infty^{(P)} \neq 0$. Let $\hat{p}$ be the model selection procedure introduced in Section 2 with $\mathcal{O} = P-1$, $c_P = c$ and $c_\mathcal{O} = 0$. Let $\tilde{\theta}$ be the corresponding post-model-selection estimator and let $G_{n,\theta,\sigma}(t|p)$ be as defined in Section 2.1. Condition (24) now implies the following: For every $\theta \in M_{P-1}$ which has exactly $P-1$ nonzero coordinates,

$$
\begin{aligned}
&\lim_{n\to\infty} P_{n,\theta,\sigma}(\{\hat{\mathfrak{r}} = \mathfrak{r}_{\text{full}}\} \blacktriangle \{\hat{p} = P\}) \\
&\quad = \lim_{n\to\infty} P_{n,\theta,\sigma}(\{\hat{\mathfrak{r}} = \mathfrak{r}_*\} \blacktriangle \{\hat{p} = P-1\}) = 0
\end{aligned}
\tag{60}
$$

holds for every $0 < \sigma < \infty$. Since the sequences $P_{n,\vartheta^{(n)},\sigma}$ and $P_{n,\theta,\sigma}$ are contiguous for $\vartheta^{(n)}$ satisfying $\|\theta - \vartheta^{(n)}\| = O(n^{-1/2})$ as remarked at the beginning of Appendix C, it follows that condition (60) continues to hold with $P_{n,\vartheta^{(n)},\sigma}$ replacing $P_{n,\theta,\sigma}$. This implies that, for every sequence of positive real numbers $s_n$ with $s_n = O(n^{-1/2})$, for every $\sigma$, $0 < \sigma < \infty$, and for every $\theta \in M_{P-1}$ which has exactly $P-1$ nonzero coordinates,

$$
\liminf_{n\to\infty} \inf_{\substack{\vartheta\in\mathbf{R}^P \\ \|\vartheta-\theta\|<s_n}} P_{n,\vartheta,\sigma}(\hat{\mathfrak{r}} = \mathfrak{r}_{\text{full}}) = \liminf_{n\to\infty} \inf_{\substack{\vartheta\in\mathbf{R}^P \\ \|\vartheta-\theta\|<s_n}} P_{n,\vartheta,\sigma}(\hat{p} = P) > 0
\tag{61}
$$

and

$$
\liminf_{n\to\infty} \inf_{\substack{\vartheta\in\mathbf{R}^P \\ \|\vartheta-\theta\|<s_n}} P_{n,\vartheta,\sigma}(\hat{\mathfrak{r}} = \mathfrak{r}_*) = \liminf_{n\to\infty} \inf_{\substack{\vartheta\in\mathbf{R}^P \\ \|\vartheta-\theta\|<s_n}} P_{n,\vartheta,\sigma}(\hat{p} = P-1) > 0,
\tag{62}
$$

hold, the positivity following from Proposition A.2. A further consequence is that

$$
\sup_{\substack{\vartheta\in\mathbf{R}^P \\ \|\vartheta-\theta\|<s_n}} \|K_{n,\vartheta,\sigma}(\cdot|\mathfrak{r}_{\text{full}}) - G_{n,\vartheta,\sigma}(\cdot|P)\|_{\text{TV}} \to 0
\tag{63}
$$

and

$$
\sup_{\substack{\vartheta\in\mathbf{R}^P \\ \|\vartheta-\theta\|<s_n}} \|K_{n,\vartheta,\sigma}(\cdot|\mathfrak{r}_*) - G_{n,\vartheta,\sigma}(\cdot|P-1)\|_{\text{TV}} \to 0
\tag{64}
$$

as $n \to \infty$. From (63)–(64), we conclude that the limit of $K_{n,\theta+\gamma/\sqrt{n},\sigma}(\cdot|\mathfrak{r}_{\text{full}})$ (with respect to total variation distance) exists and coincides with $G_{\infty,\theta,\sigma,\gamma}(\cdot|P)$. Similarly, the limit of $K_{n,\theta+\gamma/\sqrt{n},\sigma}(\cdot|\mathfrak{r}_*)$ is $G_{\infty,\theta,\sigma,\gamma}(\cdot|P-1)$. Because of (61)

and (62), we may assume that all relevant probabilities are positive (at least from a certain $n_0$ onward). Repeating the proof of Theorem 4.1 with $p = P$ and where $K_{n,\vartheta,\sigma}(t|\mathfrak{r}_{\text{full}})$ replaces $G_{n,\vartheta,\sigma}(t|P)$, as well as repeating the proof of Theorem 4.2 with $p = P - 1 = \mathcal{O}$, $q^* = P$ and where $K_{n,\vartheta,\sigma}(t|\mathfrak{r}_*)$ replaces $G_{n,\vartheta,\sigma}(t|P-1)$, gives the desired result. $\square$

PROOF OF THEOREM 3.1. Observe that (60)–(64) again hold after rearranging coordinates as in the previous proof and that

$$\lim_{n\to\infty} P_{n,\theta,\sigma}(\hat{\mathfrak{t}} = \mathfrak{r}_{\text{full}}) = \lim_{n\to\infty} P_{n,\theta,\sigma}(\hat{p} = P) > 0,$$
$$\lim_{n\to\infty} P_{n,\theta,\sigma}(\hat{\mathfrak{t}} = \mathfrak{r}_*) = \lim_{n\to\infty} P_{n,\theta,\sigma}(\hat{p} = P - 1) > 0.$$

Repeating the proof of Theorem 2.3 with $q^* = P$, with $K_{n,\vartheta,\sigma}(t|\hat{\mathfrak{t}})$ replacing $G_{n,\vartheta,\sigma}(t|\hat{p})$, and using Theorem 4.5(b) instead of Theorem 4.1(b) give the desired result. $\square$

DEPARTMENT OF STATISTICS  
YALE UNIVERSITY  
24 HILLHOUSE AVENUE  
NEW HAVEN, CONNECTICUT 06511  
USA  
E-MAIL: hannes.leeb@yale.edu

DEPARTMENT OF STATISTICS  
UNIVERSITY OF VIENNA  
UNIVERSITÄTSSTRASSE 5  
A-1010 VIENNA  
AUSTRIA  
E-MAIL: benedikt.poetscher@univie.ac.at